\documentclass[12pt]{amsart}
\usepackage{fullpage}
\usepackage{fancyhdr}
\usepackage{latexsym,amsmath,amssymb,textcomp}
\usepackage{fancyhdr}
\usepackage{geometry}
\everymath{\displaystyle}
\newtheorem{theorem}{Theorem}[section]
\newtheorem{lemma}[theorem]{Lemma}
\newtheorem{proposition}[theorem]{Proposition}
\newtheorem{corollary}[theorem]{Corollary}
\newtheorem{definition}[theorem]{Definition\rm}
\newtheorem{remark}{\it \bf Remark\/}

\newtheorem{construction}[theorem]{Construction}

\newcommand{\mylabel}[1]{\label{#1}}

\newcommand{\aatop}[2]{\genfrac{}{}{0pt}{}{#1}{#2}}

\begin{document}

\author[C. Taher]{Chadi Hassan Taher}
\address{Laboratoire J. A. Dieudonn\'e, UMR 6621
\\ Universit\'e de Nice-Sophia Antipolis\\
06108 Nice, Cedex 2, France}
\email{chaditaher@hotmail.fr}

\title[Parabolic Chern character]{Calculating the parabolic Chern character of a locally abelian parabolic bundle}

\subjclass[2000]{Primary 14F05; Secondary 19L10, 14C17}

\keywords{Parabolic bundle, Chern character, Riemann-Roch}

\begin{abstract}
We calculate the parabolic Chern character of a bundle with locally abelian
parabolic structure on a smooth strict normal crossings
divisor, using the definition in terms of Deligne-Mumford stacks. 
We obtain explicit formulas for $ch_1$, $ch_2$ and $ch_3$, and verify that these correspond to the
formulas given by Borne for $ch_1$ and Mochizuki for $ch_2$. 
\end{abstract}

\maketitle

\section{Introduction}
\mylabel{sec-intro}

Let $X$ be a smooth projective variety with a strict normal crossings divisor $D=D_1+\ldots +D_n\subset X$. 
The aim of this paper is to give  an explicit formula for the parabolic Chern character of a locally abelian parabolic bundle 
on $(X,D)$ in terms of:
\newline
---the Chern character 
of the underlying usual vector bundle, 
\newline
---the divisor components $D_{i}$ in the rational Chow groups of $X$, 
\newline
---the Chern characters
of the associated-graded pieces of the parabolic filtration along the multiple intersections of the divisor components, and
\newline
---the parabolic weights.

After giving a general formula, we  
compute explicitly the parabolic first, second, third parabolic Chern characters $ch_{1}^{Par}(E), ch_{2}^{Par}(E)$ and $ch_{3}^{Par}(E)$.

The basic idea is to use the formula given in \cite{IyerSimpson}. However, their formula did not make clear the contributions of the
different elements listed above. In order to adequatly treat this question, we start with a somewhat more general framework of 
{\em unweighted parabolic sheaves}. These are like parabolic sheaves except that the real parabolic weights are not specified.
Instead, we consider linearly ordered sets $\Sigma _i$ indexing the parabolic filtrations over the components $D_i$. 
Let $\Sigma '_i$ denote the linearly ordered set of links or adjacent pairs in $\Sigma _i$. We also call these ``risers''
as $\Sigma$ can be thought of as a set of steps. The parabolic
weights are then considered as functions $\alpha _i:\Sigma '_i \rightarrow (-1,0] \subseteq \mathbb{R}$.
This division allows us to consider separately some Chern class calculations for the unweighted structures, and then the
calculation of the parabolic Chern character using the parabolic weights.

A further difficulty stems from the fact that there are classically two different ways to give a parabolic structure:
either as a collection of sheaves included in one another; or 
by fixing a bundle $E$ (typically the zero-weight sheaf) plus a collection of filtrations of $E|_{D_i}$. The formula of
\cite{IyerSimpson} is expressed in terms of the collection of sheaves, whereas we look for a formula involving
the filtrations. Thus, our first task is to investigate the relationship between these two points of view.

An important axiom concerning the parabolic structures considered here, is that they should be {\em locally abelian}. 
This means that they should locally be direct sums of parabolic line bundles. It is a condition on the simultaneous intersection of
three or more filtrations; up to points where only two divisor components intersect, the condition is automatic. This condition has
been considered by a number of authors (Borne \cite{Borne} \cite{Borne2}, Mochizuki \cite{Mochizuki}, Iyer-Simpson \cite{IyerSimpson1},
Steer-Wren \cite{SteerWren} and others) 
and is necessary for applying the formula of \cite{IyerSimpson}.

An unweighted parabolic sheaf consists then of a collection of sheaves $E_{\sigma _1,\ldots , \sigma _n}$ with $\sigma _i\in \Sigma _i$ on $X$,
whereas an unweighted parabolic structure given by filtrations consists of a bundle $E$ on $X$ together with filtrations 
$F^i_{\sigma _i}\subset E|_{D_i}$ of the restrictions to the divisor components. In the locally abelian case, these may be related
by a long exact sequence \eqref{longsequence}:  
$$
 0\longrightarrow  E_{\sigma_{1}, \ldots ,\sigma_{n}} \rightarrow E  \rightarrow  \bigoplus_{i =
1}^{n}(\xi_{i})_{\star}(L^{i}_{\sigma_{i}}) \rightarrow \underset{i<j}{\bigoplus} (\xi_{i_{j}})_{\star}(L^{i_{j}}_{\sigma_{i}, \sigma_{j}}) \rightarrow \ldots \rightarrow L_{\sigma_{1}, \ldots ,
\sigma_{n}} \longrightarrow 0. 
$$
Where $L^{i_{1}, ..., i_{q}}_{\sigma_{i_{1}},..., \sigma_{i_{q}}}$ denote the quotient sheaves supported on intersections of  the divisors $D_{i_{1}}\cap...\cap D_{i_{q}}$.

Using this long exact sequence we get a formula \eqref{equation2ch(E)} for the Chern characters of $E_{\sigma _1,\ldots , \sigma _n}$ in terms of the Chern character
of sheaves supported on intersection of the divisors $D_{i_{1}} \cap ... \cap D_{i_{q}}$ of the form:
$$
ch^{Vb}(E_{\sigma_{1},...\sigma_{n}}) = ch^{Vb}(E)  + \sum_{q=1}^{n}(-1)^{q} \sum_{i_{1}<i_{2}<...<i_{q}}ch\left(\xi_{I,{\star}}(L^{I}_{\sigma_{i_{1}},..., \sigma_{i_{q}}})\right).
$$

The notion of parabolic weight function
is then introduced, and the main work of this paper begins: we obtain the Chern characters for the $E_{\alpha _1,\ldots ,\alpha _n}$ for any $\alpha _i\in (-1,0]$;
these are then put  into the formula of \cite{IyerSimpson}, and the result is computed. This computation requires some combinatorial manipulations with
the linearly ordered sets $\Sigma _i$ notably the associated sets of risers $\Sigma '_i$ in the ordering. It yields the following formula
\eqref{bigerequation} of Theorem \ref{bigequationtheorem}:
$$
ch^{Par}(E) = ch^{Vb}(E)e^{D}  + 
$$
$$
e^{D}.\sum_{q = 1}^{n} (-1)^{q}
\sum_{i_{1}<i_{2}<...<i_{q}}
\sum_{\lambda_{i_{j}} \in \Sigma'_{i_{j}}} ch\left(\xi_{I,{\star}}(Gr^{i_{1}, ..., i_{q}}_{\lambda_{i_{1}}, ..., \lambda_{i_{q}}})\right)
\prod_{j = 1}^{q} \left[\frac{e^{D_{i_{j}}}(1 - e^{-(\alpha_{i_{j}}(\lambda_{i_{j}}) + 1)D_{i_{j}}})}{e^{D_{i_{j}}} - 1}\right].
$$
In this formula, the associated-graded sheaves corresponding to the multiple filtrations on intersections of divisor components
$D_I=D_{i_1}\cap \cdots \cap D_{i_q}$ are denoted by $Gr^{i_{1}, ..., i_{q}}_{\lambda_{i_{1}}, ..., \lambda_{i_{q}}}$.
These are sheaves on $D_I$ but are then considered as sheaves on $X$ by the inclusion $\xi_{I,{\star}}: D_{I} \hookrightarrow X$. The Chern character
$ch\left(\xi_{I,{\star}}(Gr^{i_{1}, ..., i_{q}}_{\lambda_{i_{1}}, ..., \lambda_{i_{q}}})\right)$ is the Chern character of
the coherent sheaf on $X$. This is not satisfactory, since
we want a formula involving the Chern characters of the $Gr^{i_{1}, ..., i_{q}}_{\lambda_{i_{1}}, ..., \lambda_{i_{q}}}$ on 
$D_I$. Therefore in \S \ref{rrt} we use the Grothendieck-Riemann-Roch theorem to interchange $ch$ and $\xi_{I,{\star}}$, leading to the introduction of
Todd classes of the normal bundles of the $D_I$. 
Another difficulty is the factor of $e^D$ multiplying the term $ch^{Vb}(E)$; we would like to consider the parabolic Chern 
class as a perturbation of the Chern class of the usual vector bundle $ch^{Vb}(E)$. Using the same formula for the case of trivial 
parabolic weights, which must give back $ch^{Vb}(E)$ as an answer, allows us to rewrite the difference between $ch^{Vb}(E)$
and $ch^{Vb}(E)e^D$ in a way compatible with the rest of the formula. 
After these manipulations the formula becomes \eqref{thhequation} of Theorem \ref{maintheorem}:
$$
ch^{Par}(E) = ch^{Vb}(E) \  - 
$$
$$
e^{D}.\sum_{q = 1}^{n} (-1)^{q} \sum_{i_{1}<i_{2}<...<i_{q}}
\sum_{\lambda_{i_{j}} \in \Sigma'_{i_{j}}} \prod_{j = 1}^{q} \left(\frac{1 - e^{-D_{i_{j}}}}{D_{i_{j}}}\right).\xi_{I,{\star}}\left(ch(Gr^{i_{1}, ...,i_{q}}_{\lambda_{i_{1}}, ...,\lambda_{i_{q}}})\right) +
$$
$$
e^{D}.\sum_{q = 1}^{n} (-1)^{q} \sum_{i_{1}<i_{2}<...<i_{q}}
\sum_{\lambda_{i_{j}} \in \Sigma'_{i_{j}}} \prod_{j = 1}^{q} \left[\left(\frac{1 - e^{-\left(\alpha_{i_{j}}(\lambda_{i_{j}}) + 1\right)D_{i_{j}}}}{D_{i_{j}}}\right)\right].\xi_{I,{\star}}\left(ch(Gr^{i_{1}, ..., i_{q}}_{\lambda_{i_{1}}, ..., \lambda_{i_{q}}})\right) .
$$

Finally, we would like to compute explicitly the terms $ch_{1}^{Par}(E)$, $ch_{2}^{Par}(E)$ and $ch_{3}^{Par}(E)$. 
For these, we expand the different terms  
$$
e^D, \;\;\; 
\prod_{j = 1}^{q} \left[\left(\frac{1 - e^{-\left(\alpha_{i_{j}}(\lambda_{i_{j}}) + 1\right)D_{i_{j}}}}{D_{i_{j}}}\right)\right]  ,
\;\;\; 
\prod_{j = 1}^{q} \left(\frac{1 - e^{-D_{i_{j}}}}{D_{i_{j}}}\right)
$$
in low-degree  monomials of $D_{i_j}$, and then expand the whole formula dividing the terms up according to codimension. 
Denoting by $\mathcal{S} := \{ 1,\ldots , n\}$ the set of indices for divisor components, 
we get the following formulae:
$$
$$
$\centerdot  \ ch_{0}^{Par}(E) := rank(E).[X]$
$$
$$
$\centerdot \  ch^{Par}_{1}(E) := ch_{1}^{Vb}(E) \ - \ \sum_{i_{1} \in \mathcal{S}}\sum_{\lambda_{i_{1}} \in \sum'_{i_{1}}}\alpha_{i_{1}}(\lambda_{i_{1}}).rank(Gr^{i_{1}}_{\lambda_{i_{1}}}).[D_{i_{1}}]$
$$
$$
$\centerdot \ ch^{Par}_{2}(E) := \  ch^{Vb}_{2}(E) \ -\  \sum_{i_{1} \in \mathcal{S}}\sum_{\lambda_{i_{1}} \in \sum'_{i_{1}}} \alpha_{i_{1}}(\lambda_{i_{1}}).(\xi_{i_{1}})_{\star}\left(c_{1}^{D_{i_{1}}}(Gr^{i_{1}}_{\lambda_{i_{1}}})\right) $
$$
$$
$\hspace{2.2cm} +  \ \frac{1}{2} \ \sum_{i_{1} \in \mathcal{S}}\sum_{\lambda_{i_{1}} \in \sum'_{i_{1}}} \alpha^{2}_{i_{1}}(\lambda_{i_{1}}).rank(Gr^{i_{1}}_{\lambda_{i_{1}}}).[D_{i_{1}}]^{2}$
$$
$$
$\hspace{2.2cm} +  \sum_{i_{1} <i_{2}}
\sum_{\aatop{\lambda_{i_{1}} \in \Sigma '_{i_1}}{\lambda_{i_{2}}\in \Sigma ' _{i_2}}}  \sum_{p \in Irr(D_{i_{1}} \cap D_{i_{2}})} \alpha_{i_{1}}(\lambda_{i_{1}}).\alpha_{i_{2}}(\lambda_{i_{2}}).rank_{p}(Gr^{i_{1}, i_{2}}_{\lambda_{i_{1}}, \lambda_{i_{2}}}).[D_{p}].$

For $ch_{3}^{Par}(E)$, see Section  5.
$$
$$
 The formula for $ch^{Par}_{1}(E)$ is well-known (Seshadri {\em et al}) and, in terms of the definition of Chern classes using 
Deligne-Mumford stacks, it was shown by Borne in \cite{Borne}. The formula for 
$ch^{Par}_{2}(E)$ was given by Mochizuki in \cite{Mochizuki}, and also stated as a definition by Panov \cite{Panov}.
In both cases these coincide with our result (see the discussion on page \pageref{panovpage}). 
As far as we know, no similar formula for $ch_{3}^{Par}(E)$ has appeared in the literature.

Mochizuki defines the
Chern classes using the curvature of an adapted metric and obtains his formula as a result of a difficult curvature calculation.
It should be noted that our formula concerns the classes defined via Deligne-Mumford
stacks in the rational 
Chow groups of $X$ whereas Mochizuki's definition involving curvature can only define a class in cohomology. The identity of the two formulas shows
that the curvature definition and the stack definition give the same result up to degree $2$. 
Of course they must give the same result in general: to  prove this for the higher Chern classes this is an interesting question for further study. 

A lot of thanks to my professor Carlos Simpson for his help with this work.

\section{Unweighted parabolic structures}
\mylabel{sec-unweighted}

\subsection{Index sets}
\mylabel{sec-index}
Let $X$ be a smooth projective variety over an algebraically closed
field of characteristic zero  and let $D$ be a strict normal crossings divisor
on $X$. Write $D = D_{1} + ... + D_{n}$ where $D_{i}$ are the irreducible 
smooth components, meeting transversally. We sometimes denote by $\mathcal{S}:= \{ 1,\ldots , n\}$
the set of indices for components of the divisor $D$.

%%%%%%%%%%%%%%%%%%%%%%%%%%%%%%%%%%%%%%%%%%%%%%%%%%%%%%%%%%%%%%%%%%%%%%%%%%
\begin{definition}
\mylabel{sigma}
For $i = 1, ..., n$, 
let $\Sigma_{i}$ be finite linear ordered sets
with notations $ \eta_{i} \leq ... \leq \sigma \leq \sigma ' \leq \sigma ''
\leq ... \leq \tau_{i}$ where $\eta_{i}$ is  the smallest element of
$\Sigma_{i}$ and $\tau_{i}$ the gratest element of $\Sigma_{i}$.

Let $\Sigma '_{i}$ be the connection between the $\sigma$'s i.e
$$
\Sigma '_{i} = \{(\sigma, \sigma '),\ s.t\  \sigma < \sigma '\  and\
there \ exist \  no \ \sigma ''\  with\  \sigma < \sigma ''< \sigma'
\}.
$$ 
Consider the {\em tread functions} $ m_{+} :\Sigma'_i\rightarrow
\Sigma_i$ and $m_{-} : \Sigma'_i \rightarrow \Sigma_i$  if $\lambda =
(\sigma, \sigma') \in \Sigma'_{i} $ then $ \sigma = m_{-}(\lambda), \sigma' =
m_{+}(\lambda)$.
In the other direction, consider the {\em riser functions} $C_{+} : \Sigma _i- \{\tau _i\}
\rightarrow \Sigma'_i $ and $C_{-} : \Sigma_i - \{\eta_i\} \rightarrow
\Sigma '_i $ such that $C_{+}(\sigma) = (\sigma, \sigma ')$ where $
\sigma ' > \sigma$ the next element and $C_{-}(\sigma) = (\sigma '',
\sigma)$ where $\sigma '' < \sigma$ the next smaller element.

\end{definition}
%%%%%%%%%%%%%%%%%%%%%%%%%%%%%%%%%%%%%%%%%%%%%%%%%%%%%%%%%%%%%%%%%%%%%%%%%%%%%

One can think of the elements of $\Sigma _i$ as the steps or ``treads'' of a staircase, with $\eta_i$ and $\tau_i$ the
lower and upper landings; then $\Sigma '_i$ is the set of risers between stairs. The tread function sends a 
riser to the upper and lower treads, while the riser functions send a tread to the upper and lower risers.
The upper riser of $\tau_i$ and the lower riser of $\eta_i$ are undefined.

\subsection{Two approaches}
\mylabel{approaches}

\begin{definition}
\mylabel{uwps}
An \emph{unweighted parabolic sheaf} $E_{\cdot}$ on $(X, D)$ is a
collection of sheaves $E_{\sigma}$ indexed by multi-indices
$\sigma = (\sigma_{1}, ..., \sigma_{n})$ with $\sigma_{i} \in$
$\Sigma_{i}$, together with inclusions of sheaves
$$
E_{\sigma} \hookrightarrow E_{\sigma '}
$$
whenever $\sigma_{i} \leq  \sigma '_{i}$ (a condition which we write
as $\sigma \leq  \sigma'$ in what follows), subject to the following
hypothesis:
\begin{equation}
\label{parsheafhyp}
E_{\sigma_{1}, ..., \ \sigma_{i-1},\  \eta_{i},\  \sigma_{i+1}, ...,
\sigma_{n}} \ = \ E_{\sigma_{1}, ...,\  \sigma_{i-1},\  \tau_{i}, \
\sigma_{i+1}, ..., \ \sigma_{n}}(-D_i).
\end{equation}
\end{definition}
%%%%%%%%%%%%%%%%%%%%%%%%%%%%%%%%%%%%%%%%%%%%%%%%%%%%%%%%%%%%%%%%%%%%%%%%%%%%%%%%%%%%%%%%%%%%%%%%%%

%%%%%%%%%%%%%%%%%%%%%%%%%%%%%%%%%%%%%%%%%%%%%%%%%%%%%%%%%%%%%%%%%%%%%%%%%%%%%%%%%%%%%%%%%%%%%%%%%
\begin{construction}
\mylabel{sheaftofilt}
\end{construction}

We have inclusions of sheaves

$$
E_{\sigma_{1}, ..., \sigma_{n}}
\hookrightarrow E_{\tau_{1}, ..., \tau _{i-1},\sigma_{i}, \tau _{i+1},..., \tau_{n}}
\hookrightarrow E_{\tau_{1}, ..., \tau_{n}}.
$$
Consider the exact sequence 
$$ 
0 \longrightarrow E_{\tau_{1},...,
\tau _{i-1},\sigma_{i}, \tau _{i+1}, ..., \tau_{n}}  \longrightarrow E_{\tau_{1}, ...,
\tau_{n}} \longrightarrow E_{\tau_{1}, ..., \tau_{n}}/E_{\tau_{1},
..., \sigma_{i}, ..., \tau_{n}} \longrightarrow 0,
$$
and put \ $F_{\sigma_{i}}^{i} = \frac{E_{\tau_{1}, ..., \sigma_{i}, ...,
\tau_{n}}}{E_{\tau_{1}, ..., \eta_{i}, ..., \tau_{n}}} \subset
\frac{E_{\tau_{1}, ..., \tau_{n}}}{E_{\tau_{1}, ..., \eta_{i}, ...,
\tau_{n}}} = E\mid_{D_{i}}$ \ \ then we get the exact sequence
\begin{equation}
\label{sequencestar}
0 \longrightarrow E_{\sigma_{1}, ...,
\sigma_{n}} \longrightarrow E_{\tau_{1}, ..., \tau_{n}}
\longrightarrow \underset{i}{\bigoplus} E_{\tau_{1}, ...,
\tau_{n}}/E_{\tau_{1}, ..., \sigma_{i}, ..., \tau_{n}} =
\underset{i}{\bigoplus}\frac{ E_{\tau_{1}, ...,
\tau_{n}}/E_{\tau_{1}, ..., \eta_{i}, ...,
\tau_{n}}}{E_{\tau_{1},..., \sigma_{i}, ..., \tau_{n}}/E_{\tau_{1},
..., \eta_{i}, ..., \tau_{n}}} 
\end{equation}
which can be written as
$$
0 \longrightarrow E_{\sigma_{1}, ..., \sigma_{n}}
\longrightarrow E_{\tau_{1}, ..., \tau_{n}} \longrightarrow \
\underset{i}{\bigoplus} E\mid_{D_{i}}/F_{\sigma_{i}}^{i} .
$$
Here, to $E_{\cdot}$ we associate the usual vector bundle $E:=  E_{\tau_{1}, ..., \tau_{n}}$. 

\begin{definition}
\mylabel{psgf}
Let $E$ be  a locally free sheaf on $X$ suppose we have a filtration
denoted by $F^{i} = \{ F^{i}_{\sigma}\subseteq E\mid_{D_{i}}, \sigma
\in \Sigma_{i} \}$ of $E\mid_{D_{i}}$ where $F^{i}_{\eta_{i}} = 0$
and $F^{i}_{\tau_{i}} = E\mid_{D_{i}}$ with the remaining terms being saturated subsheaves
$$ 
0 = F^{i}_{\eta_{i}} \subseteq F^{i}_{\sigma} \subseteq ... \subseteq
F^{i}_{\tau_{i}} = E\mid_{D_{i}}
$$
for each $i = 1,..., n$. We call this a \emph{parabolic structure given by
filtrations}.
\end{definition}

The construction \ref{sheaftofilt} allows us to pass from an unweighted parabolic sheaf, to a 
parabolic structure given by filtrations. 
Suppose we are given an unweighted parabolic sheaf $E_{\sigma}$
when all the component sheaves $E_{\sigma_{1}, ..., \sigma_{n}}$ are
vector bundles.  Set $E = E_{\tau_{1}, ..., \tau_{n}}$  and
$$
E\mid_{D_{i}} = E_{\tau_{1}, ...,
\tau_{n}}/E_{\tau_{1}, ..., \tau_{i-1}, \eta_{i}, ..., \tau_{n}}.
$$

The image of $E_{\tau_{1}, ..., \tau_{i-1}, \sigma_{i}, ...,
\tau_{n}}$ in $E\mid_{D_{i}}$ is a subsheaf $F_{\sigma_{i}}^{i}$,
and we assume that it is a saturated subbundle. This gives a
parabolic structure given by filtrations.

We can also go in the opposite direction. 

\begin{construction}
\mylabel{filttosheaf}
\end{construction}

Suppose  $(E,\{ F^{i}_{\sigma} \} )$ is a parabolic structure given by filtrations. 
Consider the kernel sheaves
$$
0 \longrightarrow E_{\tau_{1},..., \sigma_{i},
..., \tau_{n}}  \longrightarrow E_{\tau_{1}, ..., \tau_{n}}
\longrightarrow E\mid_{D_{i}}/F^{i}_{\sigma_{i}},
$$
define a collection of sheaves
$$
E_{\sigma_{1}, ..., \sigma_{n}} = \bigcap _i (E_{\tau_{1}, ...,
\sigma_{i}, ..., \tau_{n}}) \subset E_{\tau_{1}, ..., \tau_{n}}
$$
with has the property that
$$
E_{\sigma _{1}, ..., \sigma _{i-1}, \eta_{i},\sigma _{i+1}, ..., \sigma _{n}} = E_{\sigma _{1}, ..., \sigma _{i-1}, \tau_{i},\sigma _{i+1}, ..., \sigma _{n}}(-D_{i}).
$$
Thus we get an unweighted parabolic sheaf.

%%%%%%%%%%%%%%%%%%%%%%%%%%%%%%%%%%%%%%%%%%%%%%%%%%%%%%%%%%%%%%%%%%%%%%%%%%%%%%%%

%%%%%%%%%%%%%%%%%%%%%%%%%%%%%%%%%%%%%%%%%%%%%%%%%%%%%%%%%%%%%%%%%%%%%%%%%%%%%%%

\subsection{Locally abelian condition.}  
An unweighted parabolic line bundle is an unweighted parabolic sheaf $F$ such that
all the $F_{\sigma}$ are line bundles. An important class of examples is obtained as follows: if $\sigma '$ is a multiindex consisting of
$\sigma '_i \in \Sigma '_i$ then we can
define an unweighted parabolic line bundle denoted 
\begin{equation}
\label{Osigma}
F := \mathcal{O}_{X}(\sigma ')
\end{equation}
by setting
$$
F_{\sigma_1, \ldots , \sigma _n} := \mathcal{O}_{X}(\sum_{i = 1}^{n}\gamma_{i}D_{i})
$$
where each $\gamma_{i}$ is equal to $-1$ or $0$; with
$\gamma _i = -1$ when $\sigma _i < \sigma '_i$ and $\gamma _i=0$ when $\sigma _i>\sigma '_i$.
Note here that the relations $<,>$ are defined between treads $\sigma _i$ and risers $\sigma '_i$.

On the other hand, if $E$ is a locally -free sheaf on $X$ then it may be considered as an unweighted parabolic sheaf (we say with trivial parabolic
 structure) by setting $E_{\sigma}$ to be $E(\sum_{i = 1}^{n}\gamma_{i}D_{i})$ for $\gamma_{i}=0$ if $\sigma _i = \tau _i$
 and $\gamma _i=-1$ otherwise. 
$$
$$
\begin{definition}
Suppose $E$ is a vector bundle on $X$ and $\sigma '$ is a multiindex consisting of $\sigma '_i \in \Sigma '_i$, we can define the unweighted parabolic bundle as follows:
$$
E(\sigma ') := E \otimes \mathcal{O}_{X}(\sigma ').
$$
\end{definition}
%%%%%%%%%%%%%%%%%%%%%%%%%%%%%%%%%%%%%%%%%%%%%%%%%%%%%%%%%%%%%%%%%%%%%%%%%%
\begin{lemma}
Any unweighted parabolic line bundle has the form $L(\sigma )$ for some $\sigma '$ and $L$ a line bundle on $X$.
\end{lemma}
%%%%%%%%%%%%%%%%%%%%%%%%%%%%%%%%%%%%%%%%%%%%%%%%%%%%%%%%%%%%%%%%%%%%%%%%%%%%%%%%%%%%%%%%%%%%%%%%%%
\begin{definition}
\mylabel{locab}
An unweighted parabolic sheaf $E_{\cdot}$, or  unweighted parabolic structure given by filtrations $(E,F^{\cdot}_{\cdot})$,  
is called  \textbf{locally abelian parabolic bundle} if in a Zariski
neighbourhood of any point $x \in X$ there is an isomorphism between
$F$ and a direct sum of  unweighted parabolic line bundles.
\end{definition}
%%%%%%%%%%%%%%%%%%%%%%%%%%%%%%%%%%%%%%%%%%%%%%%%%%%%%%%%%%%%%%%%%%%%%%%%%%%%%%%%%%

\begin{lemma}
\mylabel{locabprop}
Suppose $E_{\sigma_{1}, ..., \sigma_{n}}$ define a locally abelian
parabolic bundle  on X with respect to $(D_{1}, ..., D_{n})$. Let $E
= E_{\tau_{1}, ..., \tau_{n}}$, which is a sheaf on X. Then
$E_{\sigma}$ comes from the construction \ref{filttosheaf} as above using unique
filtrations $F^{i}_{\sigma_{i}}$ of $E\mid_{D_{i}}$ and we have the
follwing properties:

1) the $E_{\sigma_{1}, ..., \sigma_{n}}$ are locally free;

2) for each q and collection of indices $(i_{1}, ..., i_{q})$ at
each point in  the q-fold intersection $p \in D_{i_{1}}\cap...\cap
D_{i_{q}}$ the filtrations $F^{i_{1}}, ..., F^{i_{q}}$ of $E_{p}$ admit a
common splitting, hence the associated-graded
$$Gr^{F^{i_{1}}}_{j_{1}} ... Gr^{F^{i_{q}}}_{j_{q}}(E_{p})$$

is independent of the order in which it is taken;

3) the functions

$$P \mapsto {\rm rk} \   Gr^{F^{i_{1}}}_{j_{1}}
 ... Gr^{F^{i_{q}}}_{j_{q}}(E_{p})$$

are locally constant functions of $P$ on the multiple intersections
$D_{i_{1}} \cap ... \cap D_{i_{q}}$.

\end{lemma}
%%%%%%%%%%%%%%%%%%%%%%%%%%%%%%%%%%%%%%%%%%%%%%%%%%%%%%%%%%%%%%%%%%%%%

Borne \cite{Borne} shows:

%%%%%%%%%%%%%%%%%%%%%%%%%%%%%%%%%%%%%%%%%%%%%%%%%%%%%%%%%%%%%%%%%%
\begin{theorem}
\mylabel{locabchar}
Suppose given a parabolic structure which is a collection of sheaves
$E_{\sigma_{1}, ..., \sigma_{n}}$ obtained from filtrations on
bundle $E$ as above. If the sheaves satisfy condition 1), or if the filtrations satisfy 2) and 3) of the
previous lemma, then the parabolic structure is a locally abelian
parabolic bundle on (X, D).

\end{theorem}
%%%%%%%%%%%%%%%%%%%%%%%%%%%%%%%%%%%%%%%%%%%%%%%%%%%%%%%%%%%%%%%%%

Now we have two directions:
$$
$$

1-  If we have a subsheaf structure

$$  
E(-D) = E_{\eta_{1}, ..., \eta_{n}} \subseteq E_{\sigma_{1}, ...,
\sigma_{n}} \subseteq E_{\tau_{1}, ..., \tau_{n}} = E 
$$ 
we can
define the filtration structures 
$$ 
F^{i}_{\sigma_{i}} \subseteq
E\mid_{D_{i}} 
$$ 
as in Construction \ref{sheaftofilt} by using the exact sequence \eqref{sequencestar}.  Then, to
calculate the Chern character of $E_{\sigma_{1}, ..., \sigma_{n}}$
in terms of filtration structure,  we must be find an
extension for the left exact sequence \eqref{sequencestar} to a long exact sequence.
$$
$$
2- Vice versa.
$$
$$
Suppose we have a locally abelian parabolic structure $\{ F^i\}$
given by filtrations, on a vector bundle $E$ with filtrations
$$ 
0 = F^{i}_{\eta_{i}} \subseteq F^{i}_{\sigma} \subseteq ... \subseteq
F^{i}_{\tau_{i}} = E\mid_{D_{i}}.
$$

Then for $\eta_{i} \leq \sigma_{i} \leq \tau_{i}$ define the
quotient sheaves supported on $D_{i}$
$$
L^{i}_{\sigma_{i}} := \frac{E\mid_{D_{i}}}{F^{i}_{\sigma_{i}}}
$$

and the parabolic structure $E_{\sigma}$ is given by

\begin{equation}
\label{equation}
E_{\sigma_{1}, ..., \sigma_{n}} = Ker(E \longrightarrow \bigoplus_{i =
1}^{n}L^{i}_{\sigma_{i}} ).
\end{equation}
$$
$$
More generally define a family of multi-index quotient sheaves by

$$
L^{i}_{\sigma_{i}} := \frac{E\mid_{D_{i}}}{F^{i}_{\sigma_{i}}} \ on \ D_{i}
$$

$$
L^{i}_{\sigma_{i},\sigma_{j}} := \frac{E\mid_{D_{i} \cap D_{j}}}{F^{i}_{\sigma_{i}}\vert_{D_{i_{j}}} + F^{j}_{\sigma_{j}}\vert_{D_{i_{j}}}} \ on \ D_{i_{j}} = D_{i} \cap D_{j}
$$
$$.$$
$$.$$
$$
L_{\sigma_{1}, ..., \sigma_{n}} := \frac{E\mid_{D_{1} \cap ... \cap D_{n}}}{F^{1}_{\sigma_{1}} +  .... + F^{n}_{\sigma_{n}}} \ on \ D_{1} \cap D_{2}\cap
... \cap D_{n} .
$$

In these notations we have $\eta_{i} \leq \sigma_{i} \leq \tau_{i}$,  \ 
for $i = 1, ..., n$.
$$
$$
If we consider quotient sheaves as corresponding to linear subspaces
of the dual projective bundle associated to $E$, then the multiple
quotients above are multiple intersections of the $L^{i}_{\sigma_{i}}$. The formula \eqref{equation}  extends to a Koszul-style resolution of the component sheaves of the parabolic structure.
%%%%%%%%%%%%%%%%%%%%%%%%%%%%%%%%%%%%%%%%%%%%%%%%%%%%%%%%%%%%%%%%%%%%%%%%%%%%
%%%%%%%%%%%%%%%%%%%%%%%%%%%%%%%%%%%%%%%%%%%%%%%%%%%%%%%%%%%%%%%%%%%%%%%%%%%%%
\begin{lemma}
Suppose that the filtrations give a locally abelian parabolic
structure, in particular they satisfy the conditions of Lemma \ref{locabprop}.
Then for any $\eta_{i} \leq \sigma_{i} \leq \tau_{i}$ the following
sequence is well defined and exact over $X$:
\begin{equation}
\label{longsequence}
 0\longrightarrow  E_{\sigma_{1}, \ldots ,\sigma_{n}} \rightarrow E  \rightarrow  \bigoplus_{i =
1}^{n}(\xi_{i})_{\star}(L^{i}_{\sigma_{i}})  \rightarrow \underset{i<j}{\bigoplus} (\xi_{i,j})_{\star}(L^{i,j}_{\sigma_{i}, \sigma_{j}}) \rightarrow \ldots \rightarrow L_{\sigma_{1} \rightarrow \ldots ,
\sigma_{n}} \longrightarrow 0. 
\end{equation}
Where $E$ is a sheaf over $X$, \ $ \underset{i}{\bigoplus}(L^{i}_{\sigma_{i}})$ are a sheaves over
$D_{i}$,  $ \underset{i<j}{\bigoplus}(L^{i,j}_{\sigma_{i}, \sigma_{j}})$ are sheaves over $D_{i,j} = D_i\cap D_j$, etc.; 
$\xi_{i}$  denotes the closed immersion  $D_{i} \hookrightarrow X$, and $\hspace{1cm}$
$(\xi_{i})_{\star} : coh(D_{i}) \longrightarrow coh(X)$ denotes  the associated Gysin map.
The general term is a sum over $I = (i_{1}, ..., i_{q})$, where $L^I$ are sheaves over
$D_{I} = D_{i_{1}} \cap ... \cap D_{i_{q}}$ pushed forward by the Gysin map
$(\xi_{I})_{\star} : coh(D_{I}) \longrightarrow coh(X)$.
\end{lemma}
\begin{proof}
The proof in \cite{IyerSimpson} is modified to cover the unweighted case.
The maps in the exact sequence are obtained from the quotient structures of the terms with alternating signs like in the Cech complex.
We just have to prove exactness. This a local question. By the locally abelian condition, we may assume that $E$ with its filtrations is a direct sum
of rank one pieces. The formation of the sequence, and its exactness, are compatible with direct sums. Therefore we may assume that $E$ has rank one, and in fact $E \cong \mathcal{O}_{X}$.

In the case $\mathcal{O} _X(\sigma ')$ as
in \eqref{Osigma}, the vector bundle $E$ is the trivial bundle and 
the filtration steps are either $0$ or all of $ \mathcal{O}_{D_{i}}$. In particular, there is 
$\eta_{i} \leq \sigma '_{i} \leq \tau_{i}$ such that $F_{j}^{i} = \mathcal{O}_{D_{i}}$ for $ j \geq \sigma '_{i}$  and  $F_{j}^{i} = 0$ for $j < \sigma'_{i}$. Then 
$$
L^{i_{1}, ..., i_{q}}_{\sigma_{i_{1}}, ..., \sigma_{i_{q}}} = \mathcal{O}_{D_{i_{1}}, ..., D_{i_{q}}}
$$
if $\sigma_{i_{j}} < \sigma '_{i_{j}}$  for all $j = 1, ..., q$ and the quotient is zero otherwise. The sequence is defined for each multiindex $\sigma_{1}, ..., \sigma_{n}$. Up to reordering the coordinates which doesn't affect the proof, we  may assume that there is $p \in [0, n]$ such that
$\sigma_{i} < \sigma '_{i}$ for $i \leq p$ but $\sigma_{i} \geqslant \sigma '_{i}$ for $i > p$. In this case, the quotient is nonzero only when $i_{1}, ..., i_{q} \leqslant p$. Furthermore,
$$
E_{\sigma_{1}, ..., \sigma_{n}} = \mathcal{O}(-D_{1} - ... - D_{p}).
$$
In local coordinates, the divisors $D_{1}, ..., D_{p}$ are coordinate divisors. Everything is constant in the other coordinate directions which
we may ignore. The complex in question becomes 
$$
\mathcal{O}(-D_{1} - ... - D_{P}) \rightarrow \mathcal{O} \rightarrow \oplus_{1\leq i \leq p}\mathcal{O}_{D_{i}} \rightarrow \oplus_{1\leq i \leq j \leq  p}\mathcal{O}_{D_{i}\cap D_{j}} \rightarrow ... \rightarrow \mathcal{O}_{D_{1}\cap .. \cap D_{p}}.
$$
Etale locally, this is exactly the same as the exterior tensor product of $p$ copies of the resolution of $ \mathcal{O}_{ \mathbb{A}^{1}}(-D)$
on the affine line $ \mathbb{A}^{1}$ with divisor $D$ corresponding to the origin,
$$
\mathcal{O}_{ \mathbb{A}^{1}}(-D) \rightarrow \mathcal{O}_{ \mathbb{A}^{1}} \rightarrow \mathcal{O}_{D} \rightarrow 0.
$$
In particular, the exterior tensor product complex is exact except at the beginning where it resolves $ \mathcal{O}(-D_{1} - ... - D_{p})$ as required.
\end{proof}

Using the resolution of the above lemma we can compute the Chern character of $E_{\sigma_{1},...\sigma_{n}}$ in terms of the Chern character of sheaves supported on intersection of the divisors $D_{i_{1}} \cap, ..., D_{i_{q}}.$ This gives us
\begin{equation}
\label{equation2ch(E)}
ch^{Vb}(E_{\sigma_{1}, ..., \sigma_{n}}) = ch^{Vb}(E)  + \sum_{q=1}^{n}(-1)^{q} \sum_{i_{1}<i_{2}<...<i_{q}}ch\left(\xi_{I,{\star}}(L^{I}_{\sigma_{i_{1}}, ..., \sigma_{i_{q}}})\right).
\end{equation}
%%%%%%%%%%%%%%%%%%%%%%%%%%%%%%%%%%%%%%%%%%%%%%%%%%%%%%%%%%%%%%%%%%%%%%%%%%%%%%%%%%%%%%%%%%%%%%%%%%%%%%
%%%%%%%%%%%%%%%%%%%%%%%%%%%%%%%%%%%%%%%%%%%%%%%%%%%%%%%%%%%%%%%%%%%%%%%%%%%%%%%%%%%%%%%%%%%%%%%%%%%%%%%%%%%%%%

\begin{definition}
Let $I = (i_{1}, ..., i_{q})$ for $1 < i_{1} < ... < i_{q} \leq n$ and analyse the quotient 
$L_{\sigma_{i_{1}}, ..., \sigma_{i_{q}}}^{i_{1}, ..., i_{q}}$ along the multiple intersection $D_{i_{1}, ..., i_{q}}.$
There, the sheaf $E\vert_{D_{i_{1}, ..., i_{q}}}$ has $q$ filtrations $F_{\sigma_{i_{j}}}^{i_{j}}\vert_{D_{i_{1}, ..., D_{i_{q}}}}$
indexed by $\sigma_{i_{j}} \in \Sigma_{i_{j}}$ leading to a multiple associated-graded defined as follows. put 
$$
F_{\sigma_{i_{1}}, ..., \sigma_{i_{q}}}^{i_{1}, ..., i_{q}} := F_{\sigma_{i_{1}}}^{i_{1}} \cap ... \cap F_{\sigma_{i_{q}}}^{i_{q}} \subset
E\vert_{D_{i_{1}, ..., i_{q}}}.
$$
Where  $\sigma_{i_{j}} \in \Sigma_{i_{j}}$, $F^{i_{1}, ..., i_{q}}_{\eta_{i_{1}}, ..., \eta_{i_{q}}} = 0$ and $F^{i_{1}, ..., i_{q}}_{\tau_{i_{1}}, ..., \tau_{i_{q}}} = E\vert_{D_{i_{1}, ..., i_{q}}}$, we have $F^{i_{1}, ..., i_{q}}_{\sigma_{i_{1}}, ..., \sigma_{i_{q}}} \subset F^{i_{1}, ..., i_{q}}_{\sigma '_{i_{1}}, ..., \sigma '_{i_{q}}}$ if $\sigma_{i_{1}} \leqslant \sigma '_{i_{1}}$ and $\sigma_{i_{q}} \leqslant \sigma'_{i_{q}}$.
Then for a multiindex of risers $\lambda _{i_j}\in \Sigma '_{i_j}$, define
$$
Gr^{i_{1}, ..., i_{q}}_{\lambda_{i_{1}}, ..., \lambda_{i_{q}}} := \frac{F^{i_{1}, ..., i_{q}}_{m_{+}(\lambda_{i_{1}}), ..., m_{+}(\lambda_{i_{q}})}}{\sum_{j=1}^{q}F^{i_{1}, ..., i_{q}}_{m_{+}(\lambda_{i_{1}}), ..., m_{+}(\lambda_{i_{k}})-1, ..., m_{+}(\lambda_{i_{q}})}}
$$
where the indices in the denominator are almost all $m_{+}(\lambda_{i_{j}})$ but one $m_{+}(\lambda_{i_{k}})-1$.
If the parabolic structure is locally abelian then the filtrations admit a common spliting and we have
$$
Gr^{i_{1}, ..., i_{q}}_{\lambda_{i_{1}}, ..., \lambda_{i_{q}}} = Gr^{F^{i_{1}}}_{\lambda_{i_{1}}} Gr^{F^{i_{2}}}_{\lambda_{i_{2}}} ... 
Gr^{F^{i_{q}}}_{\lambda_{i_{q}}}(E\vert_{D_{i_{1}, ..., i_{q}}}).
$$
\end{definition}
$$
$$
%%%%%%%%%%%%%%%%%%%%%%%%%%%%%%%%%%%%%%%%%%%%%%%%%%%%%%%%%%%%%%%%%%%%%%%%%%%%

\begin{lemma}

Let U be a bundle over Y, and  $F_{1}^{p}, F_{2}^{p}, ..., F^{p}_{q}$ are the filtrations such that
$\exists$ common local  bases. Then in group of Grothendieck we have
$$
Gr_{F_{i}}Gr_{F_{j}} \  commute, \ and \ U \approx Gr_{F_{1}}Gr_{F_{2}} ...Gr_{F_{q}}(U).
$$
\end{lemma}
\begin{proof}
This may be proven by an inductive argument.
\end{proof}
%%%%%%%%%%%%%%%%%%%%%%%%%%%%%%%%%%%%%%%%%%%%%%%%%%%%%%%%%%%%%%%%%%%%%%%%%%%%%%
$$
$$
\begin{theorem}
Suppose given a locally abelian parabolic bundle. 
Locally over $D_{I}$ in the Zariski topology $\exists$ a finite set $\beta(\lambda_{i_{1}}, ..., \lambda_{i_{q}})$ such that we have a base over $E\vert_{D_{i}}$ of the form
$$
\lbrace e _{\lambda_{i_{1}}, ..., \lambda_{i_{q}}}; b \rbrace_{  \aatop{\lambda_{i_{j}} \in \Sigma'_{i_{j}}  }{  b \in \beta(\lambda_{i_{1}}, ..., \lambda_{i_{q}})}}  
$$
and for the filtrations $F^{i_{1}, ..., i_{q}}_{\sigma_{i_{1}}, ..., \sigma_{i_{q}}}$ admit a base  of the form
$$
\lbrace e_{\lambda_{i_{1}}, ..., \lambda_{i_{q}}}; b \rbrace_{\lambda_{i_{j}} < \sigma_{i_{j}} \ iff \ m_{+}(\lambda_{i_{j}}) \leq \sigma_{i_{j}}}
$$
and 
$$
\lbrace \overline{e}_{\lambda_{i_{1}}, ..., \lambda_{i_{q}}}; b \rbrace_{\lambda_{i_{j}} < \Sigma'_{i_{j}}}
$$
form a base of $Gr^{i_{1}, ..., i_{q}}_{\lambda_{i_{1}}, ..., \lambda_{i_{q}}}$.
\end{theorem}
\begin{proof}
By the locally abelian condition, locally we may assume that the parabolic bundle is a direct sum of parabolic line bundles. For these,
the bases have either zero or one elements and we can verify which are nonempty in terms of conditions on the $\sigma _i$. 
\end{proof}

%%%%%%%%%%%%%%%%%%%%%%%%%%%%%%%%%%%%%%%%%%%%%%%%%%%%%%%%%%%%%%%%%%%%%%%%%%%%%%%

\begin{corollary}
In the Grothendieck group of sheaves on $D_{i_{1}} \cap ... \cap D_{i_{q}}$, we have an equivalence
$$
F^{i_{1}, ..., i_{q}}_{\sigma_{i_{1}}, ..., \sigma_{i_{q}}} \approx \sum_{\aatop{ \aatop{\lambda_{i_{1}} < \sigma_{i_{1}}}{\cdots} }{  \lambda_{i_{q}} < \sigma_{i_{q}}}} Gr^{i_{1}, ..., i_{q}}_{\lambda_{i_{1}}, ..., \lambda_{i_{q}}}  \ \ \ \ \ \ \ \ \ \ \ \ \ \ \ \ \ \ \ \ \ \ \ \ \ \  (a)
$$
and
$$
\xi_{I,{\star}}\left(L^{i_{1}, ..., i_{q}}_{\sigma_{i_{1}},...,\sigma_{i_{q}}}\right) \approx  \xi_{I,{\star}}\left(\sum_{\aatop{\aatop{\sigma_{i_{1}} < \lambda_{i_{1}}}{\cdots} }{  \sigma_{i_{q}} < \lambda_{i_{q}}}} Gr^{i_{1}, ..., i_{q}}_{\lambda_{i_{1}}, ..., \lambda_{i_{q}}}\right) \ \ \ \ \ \ \ \ \ \ \ \ \ \ (b)
$$
\end{corollary}

Now apply the part $(b)$ of the above corollary  in equation \eqref{equation2ch(E)} we get a formula for the Chern character of $E_{\sigma_{1}, ..., \sigma_{n}}$ 
in terms of the associated graded as follows :
\begin{equation}
\label{medequation}
ch^{Vb}(E_{\sigma_{1}, ..., \sigma_{n}}) = ch^{Vb}(E) + \sum_{q = 1}^{n} (-1)^{q} \sum_{i_{1}<i_{2}<...<i_{q}} \sum_{\aatop{ \aatop{\sigma_{i_{1}} < \lambda_{i_{1}}}{\cdots}}{  \sigma_{i_{q}} < \lambda_{i_{q}}}}ch\left(\xi_{I,{\star}}(Gr^{i_{1}, ..., i_{q}}_{\lambda_{i_{1}}, ..., \lambda_{i_{q}}})\right). 
\end{equation}
%%%%%%%%%%%%%%%%%%%%%%%%%%%%%%%%%%%%%%%%%%%%%%%%%%%%%%%%%%%%%%%%%%%%%%%%%%%%%%%%%%%%%%%%%%
%%%%%%%%%%%%%%%%%%%%%%%%%%%%%%%%%%%%%%%%%%%%%%%%%%%%%%%%%%%%%%%%%%%%%%%%%%%%%%%%%%%
\section{Weighted parabolic structures}

The next step is to introduce the notion of {\em weight function}, providing a real number $\alpha (i)(\lambda _i)$ for each $\lambda _i\in \Sigma '_{i}$.
The weights naturally go with the ``risers'' of the linearly ordered sets, which is why we introduced the sets $\Sigma '_i$ above.

We prolong $\Sigma_{i}$ by adding its $\mathbb{Z}$-translates. Define
$$
\Phi_{i} = \mathbb{Z}.\Sigma_{i} := \mathbb{Z} \times \Sigma_{i}\diagup \thicksim
$$
s.t $(k, \tau _i) \sim (k + 1, \eta _i)$,
and
$$
\Phi_{i}' :=  \mathbb{Z} \times \Sigma_{i}'  .
$$
Prolong the functions 
$$
C_{-},C_{+} : \mathbb{Z}.\Sigma_{i} \longrightarrow \mathbb{Z} \times \Sigma_{i}'
$$
by setting $C_{+}(k, \tau) =(k + 1, C_+ (\eta ))$,
and $C_{-}(k , \eta) =(k - 1, C_-(\tau ))$.

For any unweighted parabolic sheaf we prolong the notation of 
$E_{\sigma_{1}, ..., \sigma_{n}}$  to  sheaves
$E_{\varphi_{1}, ..., \varphi_{n}}$ defined for all 
$\varphi_{i} = (k_{i}, \sigma_{i}) \in \mathbb{Z} . \Sigma_{i}$ as follows:
define 
$$
E_{\varphi_{1}, ..., \varphi_{n}} := E_{\sigma_{1}, ..., \sigma_{n}}(\Sigma k_{i}D_{i}).
$$
This is well defined modulo the equivalence relation defining $\Phi _i$, because of the
condition \eqref{parsheafhyp}.
This gives the property
$$
E_{\varphi_{1} + l_{1}, ..., \varphi_{n} + l_{n}} = E_{\varphi_{1}, ..., \varphi_{n}}(l_{1}D_{1} + ... + l_{n}D_{n})
$$
where $l_{i} \in \mathbb{Z}.$

%%%%%%%%%%%%%%%%%%%%%%%%%%%%%%%%%%%%%%%%%%%%%%%%%%%%%%%%%%%%%%%%%%%%%%%%%%%%%%%%%%%%%%%
\begin{definition}
\mylabel{weightfun}
A {\em weight function} is a collection of functions $$\alpha(i) : \Sigma'_{i} \rightarrow
(-1,0] \subseteq \mathbb{R}$$
which are increasing, i.e. $\alpha(i)(\lambda ') \leq \alpha(i)(\lambda )$ when 
$\lambda ' \leq \lambda$.
\end{definition}

To transform from unweighted parabolic structure $\longrightarrow$ weighted parabolic structure we must extend the funchtion $\alpha(i)$ 
to all of $\mathbb{Z}. \Sigma '_i$ by : 
$$
\alpha(i) : \Sigma_{i}' \longrightarrow ( -1, 0] \longrightarrow \alpha(i) : \Phi'_{i} \longrightarrow \mathbb{R}
$$
$\hspace{7.52cm}$ s.t  $ \alpha(i) (k, \sigma) = k + \alpha(i)(\sigma)$.
$$
$$
Now define  intervals by :
$$
Int\left(\alpha(i), \sigma_{i}\right) = \left(\alpha(i)C_{-}(\sigma_{i}), \alpha(i)C_{+}(\sigma_{i})\right],
$$
and 
$$
Int\left(\alpha(i), \varphi_{i}\right) = \left(\alpha(i)C_{-}(\varphi_{i}), \alpha(i)C_{+}(\varphi_{i})\right].
$$
$$
$$

We can now define the weighted or usual parabolic sheaf, associated to an unweighted parabolic sheaf and a weight function.
Consider the sheaves $E_{\beta_{1}, ..., \beta_{n}}$ for every $\beta_{i} \in \mathbb{R}^{n}$, \  given $\beta_{1}, ..., \beta_{n}\  
 \forall \ i
 \  \exists! \  \varphi_{i} \in \Phi'_{i}$ such that $\beta_{i} \in Int(\alpha(i), \varphi_{i})$, then we define
$$
E_{\beta_{1}, ..., \beta_{n}} := E_{\varphi_{1}, ..., \varphi_{n}}
$$
This defines a parabolic sheaf in the usual sense \cite{MaruyamaYokogawa} \cite{Mochizuki} \cite{Borne2} 
\cite{IyerSimpson}. 

\begin{theorem}
Suppose $E$ is a weighted parabolic bundle on $X$ with respect to $D_{1}, ..., D_{n}$. Then we have the following formula for the Chern character
of $E$ :
\begin{equation}
\label{int equation}
ch^{Par}(E) = \frac{\int_{\beta_{1}=0}^{1} ... \int_{\beta_{n}=0}^{1} e^{-\Sigma\beta_{i}D_{i}}ch^{Vb}(E_{\beta_{1}, ..., \beta_{n}})}{\int_{\beta_{1}=0}^{1} ... \int_{\beta_{n}=0}^{1} e^{-\Sigma\beta_{i}D_{i}}}.
\end{equation}
\end{theorem}
\begin{proof}
See \cite{IyerSimpson} (15), p. 35. 
\end{proof}

In this formula note the exponentials of real combinations of divisors are interpreted as formal polynomials.
The power series for the exponential terminates because the product structure of $CH^{>0}(X)$ is nilpotent.

If the weights are real, then we need the integrals as in the formula, and the result is
in $CH(X)\otimes _{\mathbb{Z}}\mathbb{R}$. If the weights are rational, then the integrals may be replaced by sums
as in \cite[Theorem 5.8]{IyerSimpson}. In this case the answer lies in 
$CH(X)\otimes _{\mathbb{Z}}\mathbb{Q}$. In what follows, if we were to replace the integrals by corresponding sums
the answer would come out the same (a factor in the numerator depending on the denominator of the
rational weights which are used, will cancel out with the same factor in the numerator).
In order to simplify notation we keep to the integral formula. 

Let $\varphi_{i} = (\sigma_{i} + 1)$ then
$$
ch^{Par}(E) = \sum_{\varphi_{1}...\varphi_{n} \in \Phi_{1} \times...\times \Phi_{n}} \frac{\int_{\beta_{1}
\in Int(\alpha(1), \varphi_{1}) \cap (0, 1]} ... \int_{\beta_{n}
\in Int(\alpha(n), \varphi_{n}) \cap (0, 1]}e^{-\Sigma \beta_{i}D_{i}} ch^{Vb}(E_{\varphi_{1},...,\varphi_{n}})}{\int_{0}^{1} ... \int_{0}^{1} e^{-\Sigma \beta_{i}D_{i}}}
$$
%%%%%%%%%%%%%%%%%%%%%%%%%%%%%%%%%%%%%%%%%%%%%%%%%%%%%%%%%%%%%%%%%%%%%%%%%%%%%%%%%%%%%%%%%%%%%%%%%%%%%%%%%%%%%%
$$
$$

\begin{remark}
$Int(\alpha(i), \varphi_{i}) \cap (0, 1] = \emptyset$  if  $\varphi_{i} \nsubseteq im\left(\lbrace1\rbrace \times \Sigma_{i}\right) \rightarrow \Phi_{i}$  i.e  just if $\varphi_{i} = \sigma_{i} +1$ for 
$\sigma_{i} \in \Sigma_{i}$.

\end{remark}
%%%%%%%%%%%%%%%%%%%%%%%%%%%%%%%%%%%%%%%%%%%%%%%%%%%%%%%%%%%%%%%%%%%%%%%%%%%%%%%%%%%%%%%%%%%%%%
\begin{definition}
Let $\varphi_{i} = \sigma_{i} + 1$, for $\sigma_{i} \in \Sigma_{i}$. Define domains by
$$
Dom(\alpha(i), \sigma_{i}) := Int(\alpha(i), \sigma_{i} + 1) \cap (0, 1],
$$
then
$$
ch^{Par}(E) = \sum_{\sigma_{1}...\sigma_{n} \in \Sigma_{1} \times ... \times \Sigma_{n}} \frac{\int_{Dom(\alpha(1), \sigma_{1})} ... \int_{Dom(\alpha(n), \sigma_{n})}e^{-\Sigma \beta_{i}D_{i}} ch^{Vb}(E_{\varphi_{1},...,\varphi_{n}})}{\int_{0}^{1} ... \int_{0}^{1}e^{-\Sigma \beta_{i}D_{i}}}.
$$
\end{definition}

We have 
$$
E_{\varphi_{1}, .., \varphi_{n}} = E_{\sigma_{1}, ..., \sigma_{n}}(D_{1} + ... + D_{n})
$$
then 
$$
ch^{Vb}(E_{\varphi_{1}, .., \varphi_{n}}) = ch^{Vb}(E_{\sigma_{1}, ..., \sigma_{n}}) e^{D_{1} + ... + D_{n}}
$$
therefore

$$
ch^{Par}(E) = \left(\frac{1}{\int_{0}^{1} ... \int_{0}^{1} e^{-\Sigma \beta_{i}D_{i}}}\right) \times
$$
$$
$$
$\sum_{\sigma_{1}...\sigma_{n} \in \Sigma_{1} \times ... \times \Sigma_{n}} ch^{Vb}(E_{\sigma_{1}, ..., \sigma_{n}}) .\left(\int_{\beta_{1} = \alpha_{-}(\sigma_{1}) + 1}^{\alpha_{+}(\sigma_{1}) + 1} .
\int_{\beta_{2} = \alpha_{-}(\sigma_{2}) + 1}^{\alpha_{+}(\sigma_{2}) + 1} ...  \int_{\beta_{n} = \alpha_{-}(\sigma_{n}) + 1}^{\alpha_{+}(\sigma_{n}) + 1} e^{-(\Sigma \beta_{i}D_{i}) + \Sigma D_{i}} d\beta\right)$
$$
$$
for $i = 1, ..., n$,  where
$$
$$ 

$\alpha_{+}(\sigma_{i})=\begin{cases}\alpha(i)(C_{+}(\sigma_{i}))  & \text{or,} \\0 &\text{if }\sigma_{i} = \tau_{i} \end{cases}$
$$
$$
and
$$
$$
$\alpha_{-}(\sigma_{i})=\begin{cases}\alpha(i)(C_{-}(\sigma_{i}))  & \text{or,} \\-1 &\text{if }\sigma_{i} = \eta_{i} \end{cases}$
$$
$$
so $Dom (\alpha(i), \sigma_{i} ) = (\alpha_{-}(\sigma_{i}) + 1, \alpha_{+}(\sigma_{i}) + 1]$.
$$
$$

Then
$$
ch^{Par}(E) = \left(\frac{1}{\int_{0}^{1} ... \int_{0}^{1} e^{-\Sigma \beta_{i}D_{i}}}\right) \times
$$
$$
$$
$\sum_{\sigma_{1}...\sigma_{n} \in \Sigma_{1} \times ... \times \Sigma_{n}} ch^{Vb}(E_{\sigma_{1}, ..., \sigma_{n}}).\left(\int_{\beta_{1} = \alpha_{-}(\sigma_{1}) + 1}^{\alpha_{+}(\sigma_{1}) + 1}.
\int_{\beta_{2} = \alpha_{-}(\sigma_{2}) + 1}^{\alpha_{+}(\sigma_{2}) + 1} ...  \int_{\beta_{n} = \alpha_{-}(\sigma_{n}) + 1}^{\alpha_{+}(\sigma_{n}) + 1} e^{-\Sigma (\beta_{i} - 1)D_{i}} d\beta\right) .$
$$
$$

$$
$$
Take 
$\gamma_{i} = \beta_{i} - 1$ for $i = 1, ..., n$
then
$$
ch^{Par}(E) = \left(\frac{1}{\int_{0}^{1} ... \int_{0}^{1} e^{-\Sigma \beta_{i}D_{i}}}\right) \times
$$
$\sum_{\sigma_{1}...\sigma_{n} \in \Sigma_{1} \times .. \times \Sigma_{n}} ch^{Vb}(E_{\sigma_{1}, ..., \sigma_{n}})\left(\int_{\gamma_{1} = \alpha_{-}(\sigma_{1})}^{\alpha_{+}(\sigma_{1})}e^{-\gamma_{1}D_{1}}d\gamma_{1} . \int_{\gamma_{2} = \alpha_{-}(\sigma_{2})}^{\alpha_{+}(\sigma_{2})}e^{-\gamma_{2}D_{2}}d\gamma_{2}\ ..\int_{\gamma_{n} = \alpha_{-}(\sigma_{n})}^{\alpha_{+}(\sigma_{n})}e^{-\gamma_{n}D_{n}}d\gamma_{n}\right)$
$$  
$$
and we have 
$$
ch^{Vb}(E_{\sigma_{1}, ..., \sigma_{n}}) = ch^{Vb}(E) + \sum_{q = 1}^{n} (-1)^{q} \sum_{i_{1}<i_{2}<...<i_{q}} \sum_{\aatop{\sigma_{i_{1}} < \lambda_{i_{1}} }{  \sigma_{i_{q}} < \lambda_{i_{q}}}}ch\left(\xi_{I,{\star}}(Gr^{i_{1}, ..., i_{q}}_{\lambda_{i_{1}}, ..., \lambda_{i_{q}}})\right). 
$$

Then
$$
ch^{Par}(E) = \left(\frac{1}{\int_{0}^{1} ... \int_{0}^{1} e^{-\Sigma \beta_{i}D_{i}}}\right) \times 
$$
$$
$$
$ \hspace{0.7cm}\sum_{\sigma_{1}...\sigma_{n} \in \Sigma_{1} \times ... \times \Sigma_{n}} \left(ch^{Vb}(E) + \sum_{q = 1}^{n} (-1)^{q} \sum_{i_{1}<i_{2}<...<i_{q}}\sum_{\aatop{\sigma_{i_{1}} < \lambda_{i_{1}} }{  \sigma_{i_{q}} < \lambda_{i_{q}}}}ch\left(\xi_{I,{\star}}(Gr^{i_{1}, ..., i_{q}}_{\lambda_{i_{1}}, ..., \lambda_{i_{q}}})\right)\right) \times $
$$
$$
$\hspace{2.4cm}\left(\int_{\gamma_{1} = \alpha_{-}(\sigma_{1})}^{\alpha_{+}(\sigma_{1})}e^{-\gamma_{1}D_{1}}d\gamma_{1}  \int_{\gamma_{2} = \alpha_{-}(\sigma_{2})}^{\alpha_{+}(\sigma_{2})}e^{-\gamma_{2}D_{2}}d\gamma_{2} \  ... \  \int_{\gamma_{n} = \alpha_{-}(\sigma_{n})}^{\alpha_{+}(\sigma_{n})}e^{-\gamma_{n}D_{n}}d\gamma_{n}\right)$  
$$
$$
$$
$$
$$
= \left(\frac{ch^{Vb}(E)}{\int_{0}^{1} ... \int_{0}^{1} e^{-\Sigma \beta_{i}D_{i}}}\right) \times
$$
$$
$$
$\left(\sum_{\sigma_{1}\in \Sigma _{1}}\int_{\gamma_{1} = \alpha_{-}(\sigma_{1})}^{\alpha_{+}(\sigma_{1})}e^{-\gamma_{1}D_{1}}d\gamma_{1}\right)\cdot 
\left(\sum_{\sigma_{2}\in \Sigma_{2}}\int_{\gamma_{2} = \alpha_{-}(\sigma_{2})}^{\alpha_{+}(\sigma_{2})}e^{-\gamma_{2}D_{2}}d\gamma_{2}\right)
 \cdots \left(\sum_{\sigma_{n}\in \Sigma_{n}}\int_{\gamma_{n} = \alpha_{-}(\sigma_{n})}^{\alpha_{+}(\sigma_{n})}e^{-\gamma_{n}D_{n}}d\gamma_{n}\right)$
$$
$$
$$
$$
$$
 + \ \left(\frac{1}{\int_{0}^{1} ... \int_{0}^{1} e^{-\Sigma \beta_{i}D_{i}}}\right).\left[ \sum_{q = 1}^{n} (-1)^{q} \sum_{i_{1}<i_{2}<...<i_{q}}
\sum_{\lambda_{i_{1}} \ldots  \lambda_{i_{q}}} ch\left(\xi_{I,{\star}}(Gr^{i_{1}, ..., i_{q}}_{\lambda_{i_{1}}, ..., \lambda_{i_{q}}})\right)\right. \times
$$
$$
$$
$$
 \hspace{1cm}\left(\sum_{\aatop{\sigma_{1} \in \Sigma_{1} }{ \sigma_{1} < \lambda_{1} \  if \  1 \in I}}\int_{\gamma_{1} = \alpha_{-}(\sigma_{1})}^{\alpha_{+}(\sigma_{1})}e^{-\gamma_{1}D_{1}}d\gamma_{1}\right)\cdot \left(\sum_{\aatop{\sigma_{2} \in \Sigma_{2} }{ \sigma_{2} < \lambda_{2} \  if \  2 \in I}} \int_{\gamma_{2} = \alpha_{-}(\sigma_{2})}^{\alpha_{+}(\sigma_{2})}e^{-\gamma_{2}D_{2}}d\gamma_{2}\right) \cdots 
$$
$$
$$
$$
\left. \cdot \left(\sum _{\aatop{\sigma_{n} \in \Sigma_{n} }{ \sigma_{n} < \lambda_{n} \  if \  n \in I}}\int_{\gamma_{n} = \alpha_{-}(\sigma_{n})}^{\alpha_{+}(\sigma_{n})}e^{-\gamma_{n}D_{n}}d\gamma_{n}\right) \right]
$$
$$
$$
which can be written as
$$
ch^{Par}(E) = A + B
$$
$$
$$
where
$$  
 A =  \left(\frac{ch^{Vb}(E)}{\int_{0}^{1} ... \int_{0}^{1} e^{-\Sigma \beta_{i}D_{i}}}\right) \times
$$
$$
$$
$\left(\sum_{\sigma_{1}\in \Sigma _{1}}\int_{\gamma_{1} = \alpha_{-}(\sigma_{1})}^{\alpha_{+}(\sigma_{1})}e^{-\gamma_{1}D_{1}}d\gamma_{1}\right)\cdot
\left(\sum_{\sigma_{2}\in \Sigma _{2}}\int_{\gamma_{2} = \alpha_{-}(\sigma_{2})}^{\alpha_{+}(\sigma_{2})}e^{-\gamma_{2}D_{2}}d\gamma_{2}\right)
 ... \left(\sum_{\sigma_{n}\in \Sigma _{n}}\int_{\gamma_{n} = \alpha_{-}(\sigma_{n})}^{\alpha_{+}(\sigma_{n})}e^{-\gamma_{n}D_{n}}d\gamma_{n}\right)$
$$
$$
One can note that 
$\int_{a}^{b} e^{-\rho D}d\rho = \frac{e^{-aD} - e^{-bD}}{D} = e^{-aD}\frac{(1 - e^{(a -b)D})}{D} =  e^{-aD}\frac{(a-b)}{td((a -b)D)}$
where td is the Todd class.
$$
$$
We have
$\sum_{\sigma_{i}\in \Sigma _{i}}\int_{\gamma_{i} = \alpha_{-}(\sigma_{i})}^{\alpha_{+}(\sigma_{i})}e^{-\gamma_{i}D_{i}}d\gamma_{i} =
\int_{-1}^{0}e^{-\gamma_{i}D_{i}}d\gamma_{i}$ for $i =1,2, ..., n$ and put $\beta = \gamma$ as integration variable
$$
$$
then
$$
A = \left(\frac{ch^{Vb}(E)}{\int_{0}^{1} ... \int_{0}^{1} e^{-\Sigma \gamma_{i}D_{i}}}\right).\left(\int_{-1}^{0}e^{-\gamma_{1}D_{1}}d\gamma_{1}\right).\left(\int_{-1}^{0}e^{-\gamma_{2}D_{2}}d\gamma_{2}\right) ... \left(\int_{-1}^{0}e^{-\gamma_{n}D_{n}}d\gamma_{n}\right)
$$
$$
$$
$ =\left[ \left(\frac{\int_{-1}^{0}e^{-\gamma_{1}D_{1}}d\gamma_{1}}{\int_{0}^{1}e^{-\gamma_{1}D_{1}}d\gamma_{1}}\right). \left(\frac{\int_{-1}^{0}e^{-\gamma_{2}D_{2}}d\gamma_{2}}{\int_{0}^{1}e^{-\gamma_{2}D_{2}}d\gamma_{2}}\right)... \left(\frac{\int_{-1}^{0}e^{-\gamma_{n}D_{n}}d\gamma_{n}}{\int_{0}^{1}e^{-\gamma_{n}D_{n}}d\gamma_{n}}\right)\right].ch^{Vb}(E)$
$$
$$
where 

$$
A_{i} = \frac{\int_{-1}^{0}e^{-\gamma_{i}D_{i}}d\gamma_{i}}{\int_{0}^{1}e^{-\gamma_{i}D_{i}}d\gamma_{i}} = 
\frac{e^{-D_{i}} - 1}{1 - e^{-D_{i}}} = \frac{e^{D_{i}}(1 - e^{-D_{i}})}{1 - e^{-D_{i}}} = e^{D_{i}}
$$
$$
$$
$$
for \hspace{0.5cm}   i = 1, 2, ..., n
$$
therefore
$$
A = ch^{Vb}(E).e^{D_{1}}.e^{D_{2}}...e^{D_{n}} = ch^{Vb}(E)e^{D}
$$
where $D = D_{1} + D_{2} + ... + D_{n}$
$$
$$
and
$$
B = \left(\frac{1}{\int_{0}^{1} ... \int_{0}^{1} e^{-\Sigma \beta_{i}D_{i}}}\right).\left[ \sum_{q = 1}^{n} (-1)^{q} \sum_{i_{1}<i_{2}<...<i_{q}}
\sum_{\aatop{\lambda_{i_{1}} }{  \lambda_{i_{q}}}} ch\left(\xi_{I,{\star}}(Gr^{i_{1}, ..., i_{q}}_{\lambda_{i_{1}}, ..., \lambda_{i_{q}}})\right) \times \right.  
$$
$$
$$
$$
\hspace{1cm}\left(\sum_{\aatop{\sigma_{1} \in \Sigma_{1} }{ \sigma_{1} < \lambda_{1} \ if \  1 \in I}}\int_{\gamma_{1} = \alpha_{-}(\sigma_{1})}^{\alpha_{+}(\sigma_{1})}e^{-\gamma_{1}D_{1}}d\gamma_{1}\right) \cdot
\left(\sum_{\aatop{\sigma_{2} \in \Sigma_{2} }{ \sigma_{2} < \lambda_{2} \ if \  2 \in I}} \int_{\gamma_{2} = \alpha_{-}(\sigma_{2})}^{\alpha_{+}(\sigma_{2})}e^{-\gamma_{2}D_{2}}d\gamma_{2}\right)
$$
$$
$$
$$
\left. ...  \left(\sum _{\aatop{\sigma_{n} \in \Sigma_{n} }{ \sigma_{n} < \lambda_{n} \ if \  n \in I }}\int_{\gamma_{n} = \alpha_{-}(\sigma_{n})}^{\alpha_{+}(\sigma_{n})}e^{-\gamma_{n}D_{n}}d\gamma_{n}\right) \right]
$$
$$
$$
The sums of integrals can be expressed as single integrals; 
if $i \notin I$, set formally $\alpha_{i}(\lambda_{i})_{} := 0$ in the following expression:
$$
$$

$$
B = \left[\frac{(\int_{-1}^{0}e^{-\gamma_{1}D_{1}}d\gamma_{1})(\int_{-1}^{0}e^{-\gamma_{2}D_{2}}d\gamma_{2}) ... (\int_{-1}^{0}e^{-\gamma_{n}D_{n}}d\gamma_{n})}{(\int_{0}^{1}e^{-\gamma_{1}D_{1}}d\gamma_{1})(\int_{0}^{1}e^{-\gamma_{2}D_{2}}d\gamma_{2}) ... (\int_{0}^{1}e^{-\gamma_{n}D_{n}}d\gamma_{n})}\right]\times
$$
$$
$$
$$ 
\hspace{1cm}\left( \sum_{q = 1}^{n} (-1)^{q} \sum_{i_{1}<i_{2}<...<i_{q}}
\sum_{\lambda_{i_{j}} \in \Sigma'_{i_{j}} } ch\left(\xi_{I,{\star}}(Gr^{i_{1}, ..., i_{q}}_{\lambda_{i_{1}}, ..., \lambda_{i_{q}}})\right)\right. \times
$$
$$
$$
$$
\hspace{1cm}\left[\frac{(\int_{-1}^{0}e^{-\gamma_{1}D_{1}}d\gamma_{1})(\int_{-1}^{0}e^{-\gamma_{2}D_{2}}d\gamma_{2}) ... (\int_{-1}^{0}e^{-\gamma_{n}D_{n}}d\gamma_{n})}{(\int_{-1}^{0}e^{-\gamma_{1}D_{1}}d\gamma_{1})(\int_{-1}^{0}e^{-\gamma_{2}D_{2}}d\gamma_{2}) ... (\int_{-1}^{0}e^{-\gamma_{n}D_{n}}d\gamma_{n})}\right] \times
$$
$$
$$
$$
\left. \hspace{1cm}\left[\frac{(\int_{-1}^{\alpha_{1}(\lambda_{1})}e^{-\gamma_{1}D_{1}}d\gamma_{1})(\int_{-1}^{\alpha_{2}(\lambda_{2})}e^{-\gamma_{2}D_{2}}d\gamma_{2}) ... (\int_{-1}^{\alpha_{n}(\lambda_{n})}e^{-\gamma_{n}D_{n}}d\gamma_{n})}{(\int_{-1}^{0}e^{-\gamma_{1}D_{1}}d\gamma_{1})(\int_{-1}^{0}e^{-\gamma_{2}D_{2}}d\gamma_{2}) ... (\int_{-1}^{0}e^{-\gamma_{n}D_{n}}d\gamma_{n})}\right] \right)
$$
$$
$$
$$
$$
$$
= e^{D}.\sum_{q = 1}^{n} (-1)^{q} \sum_{i_{1}<i_{2}<...<i_{q}}
\sum_{\lambda_{i_{j}} \in \Sigma'_{i_{j}}} ch\left(\xi_{I,{\star}}(Gr^{i_{1}, ..., i_{q}}_{\lambda_{i_{1}}, ..., \lambda_{i_{q}}})\right)
\prod_{j = 1}^{q}\left[\frac{\int_{-1}^{\alpha_{i_{j}(\lambda_{i_{j}})}}e^{-\rho D_{i_{j}}}d\rho}{\int_{-1}^{0}e^{-\rho D_{i_{j}}}d\rho}\right]
$$
$$
= e^{D}.\sum_{q = 1}^{n} (-1)^{q} \sum_{i_{1}<i_{2}<...<i_{q}}
\sum_{\lambda_{i_{j}} \in \Sigma'_{i_{j}}} ch\left(\xi_{I,{\star}}(Gr^{i_{1}, ..., i_{q}}_{\lambda_{i_{1}}, ..., \lambda_{i_{q}}})\right)
\prod_{j = 1}^{q} \left[\frac{e^{D_{i_{j}}} - e^{-(\alpha_{i_{j}}(\lambda_{i_{j}}))D_{i_{j}}}}{e^{D_{i_{j}}} - 1}\right]
$$
$$
= e^{D}.\sum_{q = 1}^{n} (-1)^{q} \sum_{i_{1}<i_{2}<...<i_{q}}
\sum_{\lambda_{i_{j}} \in \Sigma'_{i_{j}}} ch\left(\xi_{I,{\star}}(Gr^{i_{1}, ..., i_{q}}_{\lambda_{i_{1}}, ..., \lambda_{i_{q}}})\right)
\prod_{j = 1}^{q} \left[\frac{e^{D_{i_{j}}}(1 - e^{-(\alpha_{i_{j}}(\lambda_{i_{j}}) + 1)D_{i_{j}}})}{e^{D_{i_{j}}} - 1}\right]
$$
$$
$$
therefore we have proven the following
$$
$$
\begin{theorem}
\label{bigequationtheorem}
\begin{equation}
\label{bigerequation}
ch^{Par}(E) = ch^{Vb}(E)e^{D}  + 
$$
$$
e^{D}.\sum_{q = 1}^{n} (-1)^{q}
\sum_{i_{1}<i_{2}<...<i_{q}}
\sum_{\lambda_{i_{j}} \in \Sigma'_{i_{j}}}ch\left(\xi_{I,{\star}}(Gr^{i_{1}, ..., i_{q}}_{\lambda_{i_{1}}, ..., \lambda_{i_{q}}})\right)
\prod_{j = 1}^{q} \left[\frac{e^{D_{i_{j}}}(1 - e^{-(\alpha_{i_{j}}(\lambda_{i_{j}}) + 1)D_{i_{j}}})}{e^{D_{i_{j}}} - 1}\right].
\end{equation}
\end{theorem}
$$
$$

\subsection{Riemann-Roch theorem}
\label{rrt}
The next step is to use Riemann-Roch theory to interchange $ch$ and $\xi _{I,\star}$. Let $K^{\circ}D_{I}$ denote the Grothendieck group of vector bundles(locally free sheaves) on $D_{I}$.  Each vector bundle
$E$ determines an element, denoted $[E]$, in  $K^{\circ}D_{I}$. $K^{\circ}D_{I}$ is the free abelian group on the set of isomorphism
classes of vector bundles, modulo the relations
$$
[E] \ = \  [E']  +  [E'']
$$
whenever $E'$ is a subbundle of a vector bundle $E$ with quotient bundle $E'' = E/E'$. The tensor product
makes $K^{\circ}D_{I}$ a ring: $[E].[F] = [E \otimes F]$.
 
 %%%%%%%%%%%%%%%%%%%%%%%%%%%%%%%%%%%%%%%%%%%%%%%%%%%%%%%%%%%%%%%%%%%%%%%%%%%%%%%%%%%%%%%%%%%%%%%%%%%%%%
\begin{definition}
 For any morphism $\xi_{I} :X \rightarrow D_{I}$
 there is an induced homomorphism
 $$
 \xi^{\star}_{I} : K^{\circ}D_{I} \rightarrow K^{\circ}X,
 $$ 
 taking $[E]$ to $[\xi^{\star}_{I}E]$, where $\xi^{\star}_{I}E$ is the pull-back bundle.
 \end{definition}
 
 %%%%%%%%%%%%%%%%%%%%%%%%%%%%%%%%%%%%%%%%%%%%%%%%%%%%%%%%%%%%%%%%%%%%%%%%%%%%%%%%%%%%%%%%%%%%%%%%%%%%%%%%
 
 The Grothendieck group of coherent sheaves on $D_{I}$, denoted by $K_{\circ}D_{I}$, is defined to be the free abelian
 group on the isomorphism class $[\mathcal{F}]$ of coherent sheaves on $D_{I}$, modulo the relaions
 
$$
 [\mathcal{F}] \ = \ [\mathcal{F'}] +  [\mathcal{F''}]
 $$
 for each exact sequence
 
 $$
 0 	\rightarrow \mathcal{F'} \rightarrow \mathcal{F} \rightarrow \mathcal{F''} \rightarrow 0
 $$
 of coherent sheaves on $D_{I}$. Tensor product makes $K_{\circ}D_{I}$ a $K^{\circ}D_{I}$-module:
 
$$
K^{\circ}D_{I} \otimes K_{\circ}D_{I} \rightarrow K_{\circ}D_{I},
$$
$[E].[\mathcal{F}] = [E \otimes_{\mathcal{O_{D_{I}}}}\mathcal{F}]$.

\begin{definition}

For any proper morphism $\xi_{I} : D_{I} \rightarrow X$, there is a homomorphism 
$$
\xi_{I,{\star}} : K_{\circ}D_{I} \rightarrow K_{\circ}X
$$
which takes $[\mathcal{F}]$ to $\sum_{i \geqslant 0}(-1)^{i}[R^{i}\xi_{I,{\star}}\mathcal{F}]$.where
$R^{i}\xi_{I,{\star}}\mathcal{F}$ is Grothendieck  higher direct image sheaf, the sheaf associated to the presheaf
$$
U \rightarrow H^{i}(\xi^{-1}_{I}(U), \mathcal{F})
$$
on $X$. 

\end{definition}

It is a basic fact the $R^{i}\xi_{I,{\star}}\mathcal{F}$  are coherent when $\mathcal{F}$ is coherent and 
 $\xi_{I}$ is proper. The fact that this push-forward $\xi_{I,{\star}}$ is well-defined on $K_{\circ}D_{I}$ results from the long excat
 cohomology sequence for the $R^{i}\xi_{I,{\star}}$. 
 
\begin{proposition}
 The push-forward and pull-back are related by the usual projection formula:

 $$
 \xi_{I,{\star}}(\xi_{I}^{\star}a.b) \ = \ a.(\xi_{I,{\star}}b)
 $$ 
 $$
 $$
 for  $\xi_{I} : D_{I} \rightarrow X$ proper, $a \in K^{\circ}X$, $b \in K_{\circ}D_{I}$. 
\end{proposition}

%%%%%%%%%%%%%%%%%%%%%%%%%%%%%%%%%%%%%%%%%%%%%%%%%%%%%%%%%%%%%%%%%%%%%%%

\begin{theorem}
On any $D_{I}$  there is a canonical duality homomorphism
$$
K^{\circ}D_{I} \rightarrow K_{\circ}D_{I}
$$
which takes a vector bundle to its sheaf of sections. When $D_{I}$ is non-singular, 
this duality map is an isomorphism.

\end{theorem}
%%%%%%%%%%%%%%%%%%%%%%%%%%%%%%%%%%%%%%%%%%%%%%%%%%%%%%%%%%%%%%%%%%%%%%%%%%%%%

\begin{definition}
Consider $D_{I}$ which are smooth over a given ground field $\mathbb{C}$. For such $D_{I}$ we identify $K^{\circ}D_{I}$ and $K_{\circ}D_{I}$, and write simply $K(D_{I})$.  There is a homomorphism, called the Chern character
$$
ch: K(D_{I}) \rightarrow A(D_{I})_{\mathbb{Q}}
$$
determined by the following properties:

i) ch is a homomorphism of rings;

ii) if $\xi_{I} : X \rightarrow D_{I}$, $ch \circ \xi_{I}^{\star} =  \xi_{I}^{\star} \circ ch$;

iii) if $l$ is a line bundle on $D_{I}$,
$ch[D_{I}] = exp(c_{1}(l)) = \sum_{i\geqslant0}(1/i!)c_{1}(l)^{i}$
\end{definition}

%%%%%%%%%%%%%%%%%%%%%%%%%%%%%%%%%%%%%%%%%%%%%%%%%%%%%%%%%%%%%%%%%%%%%%%%%%%%%%
\begin{theorem} 
Let $\xi_{I} : D_{I} \hookrightarrow X$ be a  smooth projective morphism of non-singular quasi-projective varieties. Then for any 
$Gr_{\lambda}^{I} \in K(D_{I})$ we have
$$
\xi_{I,{\star}}\left(ch(Gr^{i_{1}, ..., i_{q}}_{\lambda_{i_{1}}, ..., \lambda_{i_{q}}}).td(T_{D_{I}})\right) = td(T_{X}).ch\left(\xi_{I,{\star}}(Gr^{i_{1}, ..., i_{q}}_{\lambda_{i_{1}}, ..., \lambda_{i_{q}}})\right) 
$$
in $A(X)\otimes \mathbb{Q}$. where $Td(X) = td(T_{X}) \in A(X)_{\mathbb{Q}}$ is the relative tangent sheaf of $\xi_{I}$.
\end{theorem}

\noindent
See \cite[pp 286-287]{Fulton}.
%%%%%%%%%%%%%%%%%%%%%%%%%%%%%%%%%%%%%%%%%%%%%%%%%%%%%%%%%%%%%%%%%%%%%%%%%%%%%%%%
%%%%%%%%%%%%%%%%%%%%%%%%%%%%%%%%%%%%%%%%%%%%%%%%%%%%%%%%%%%%%%%%%%%%%%%%%%%%%%%%
\begin{theorem}
If $D_{I}$ is a non-singular variety set 
$$
Td(D_{I}) = td(T_{D_{I}}) \in A(D_{I})_{\mathbb{Q}}
$$
then
$$
$$
if $\xi_{I} : D_{I} \hookrightarrow X$, is a closed imbedding of codimension $q$, and $D_{I}$ is the intersection of $q$ Cartier
divisors $D_{i_{1}}, ..., D_{i_{q}}$ on $X$, then
$$
Td(D_{I}) = \xi^{\star}_{I}\left[Td(X).\prod_{j = 1}^{q} \left(\frac{1 - e^{-Di_{j}}}{D_{i_{j}}}\right)\right]. 
$$
\end{theorem}

\noindent
See \cite[p. 293]{Fulton}.

%%%%%%%%%%%%%%%%%%%%%%%%%%%%%%%%%%%%%%%%%%%%%%%%%%%%%%%%%%%%%%%%%%%%%%%%%%%%%%%%%%
$$
$$
So

$\xi_{I,{\star}}\left(ch(Gr^{i_{1}, ..., i_{q}}_{\lambda_{i_{1}}, ..., \lambda_{i_{q}}}).td(T_{D_{I}})\right) =
\xi_{I,{\star}}\left[ ch(Gr^{i_{1}, ..., i_{q}}_{\lambda_{i_{1}}, ..., \lambda_{i_{q}}}).
\xi^{\star}_{I}\left( Td(X).\prod_{j = 1}^{q}\left(\frac{1 - e^{-Di_{j}}}{D_{i_{j}}} \right)\right)\right]$
$$
$$
$$
\hspace{5,6cm}= \ Td(X).\prod_{j = 1}^{q} \left[\frac{1 - e^{-Di_{j}}}{D_{i_{j}}}\right].\xi_{I,{\star}}\left(ch(Gr^{i_{1}, ..., i_{q}}_{\lambda_{i_{1}}, ..., \lambda_{i_{q}}})\right).
$$
$$
$$
%%%%%%%%%%%%%%%%%%%%%%%%%%%%%%%%%%%%%%%%%%%%%%%%%%%%%%%%%%%%%%%%%%%%%%%%%%%%%%%%%
\begin{corollary}
\label{RRcorollary}
$$
ch\left(\xi_{I_\star}(Gr^{i_{1}, ..., i_{q}}_{\lambda_{i_{1}}, ..., \lambda_{i_{q}}})\right)  = \prod_{j = 1}^{q} \left[\frac{1 - e^{-Di_{j}}}{D_{i_{j}}}\right].\xi_{I,{\star}}\left(ch(Gr^{i_{1}, ..., i_{q}}_{\lambda_{i_{1}}, ..., \lambda_{i_{q}}})\right)
$$
$$
$$
$\hspace{5,2cm} = \prod_{j = 1}^{q} \left[\left(\frac{e^{-D_{i_{j}}}(e^{D_{i_{j}}} - 1)} {D_{i_{j}}}\right)\right].\xi_{I,{\star}}\left(ch(Gr^{i_{1}, ..., i_{q}}_{\lambda_{i_{1}}, ..., \lambda_{i_{q}}})\right). 
$
\end{corollary}
$$
$$
%%%%%%%%%%%%%%%%%%%%%%%%%%%%%%%%%%%%%%%%%%%%%%%%%%%%%%%%%%%%%%%%%%%%%%%%%%%%%%%%%%%%%
\begin{theorem}
Apply  the above corollary in the equation \eqref{bigerequation} of Theorem \ref{bigequationtheorem} we get

\begin{equation}
\label{theoequation}
ch^{Par}(E) = ch^{Vb}(E)e^{D} +
$$
$$
e^{D}.\sum_{q = 1}^{n} (-1)^{q} \sum_{i_{1}<i_{2}<...<i_{q}}
\sum_{\lambda_{i_{j}} \in \Sigma'_{i_{j}}} \prod_{j = 1}^{q} \left[\left(\frac{1 - e^{-\left(\alpha_{i_{j}}(\lambda_{i_{j}}) + 1\right)D_{i_{j}}}}{D_{i_{j}}}\right)\right].\xi_{I,{\star}}\left(ch(Gr^{i_{1}, ..., i_{q}}_{\lambda_{i_{1}}, ..., \lambda_{i_{q}}})\right).
\end{equation}
\end{theorem}
$$
$$
Now return to the equation $(\ref{medequation})$ but apply it to the bottom value of $\sigma _i = \eta _i$.
Recall that $E_{\eta _{1}, ..., \eta _{n}} = E(-D)$; but on the other hand $\eta _{i_j}< \lambda _{i_j}$ for any
$\lambda _{i_j}\in \Sigma ' _{i_j}$. Thus we get 
$$
$$
$$
ch^{Vb}(E(-D)) =  ch^{Vb}(E) + \sum_{q = 1}^{n} (-1)^{q} \sum_{i_{1}<i_{2}<...<i_{q}} \sum_{ \lambda_{i_{j}}\in \Sigma '_{i_j} }
ch\left(\xi_{I,{\star}}(Gr^{i_{1}, ..., i_{q}}_{\lambda_{i_{1}}, ..., \lambda_{i_{q}}})\right). 
$$
$\Longrightarrow$
\begin{equation}
\label{chvbequation}
ch^{Vb}(E) =  ch^{Vb}\left(E(-D)\right) -  \sum_{q = 1}^{n} (-1)^{q} \sum_{i_{1}<i_{2}<...<i_{q}} \sum_{\lambda_{i_{j}}\in \Sigma '_{i_j}}
ch\left(\xi_{I,{\star}}(Gr^{i_{1}, ..., i_{q}}_{\lambda_{i_{1}}, ..., \lambda_{i_{q}}})\right). 
\end{equation}
$$
$$
Put equation ($ \ref{chvbequation}$) in  ($ \ref{theoequation}$)  and use Corollary \ref{RRcorollary}, to get
$$
ch^{Par}(E) = ch^{Vb}\left(E(-D)\right)e^{D} - 
$$
$$
e^{D}.\sum_{q = 1}^{n} (-1)^{q} \sum_{i_{1}<i_{2}<...<i_{q}}
\sum_{\lambda_{i_{j}} \in \Sigma'_{i_{j}}} \prod_{j = 1}^{q} \left(\frac{1 - e^{-D_{i_{j}}}}{D_{i_{j}}}\right).\xi_{I,{\star}}\left(ch(Gr^{i_{1}, ...,i_{q}}_{\lambda_{i_{1}}, ...,\lambda_{i_{q}}})\right) +
$$
$$
e^{D}.\sum_{q = 1}^{n} (-1)^{q} \sum_{i_{1}<i_{2}<...<i_{q}}
\sum_{\lambda_{i_{j}} \in \Sigma'_{i_{j}}} \prod_{j = 1}^{q} \left[\left(\frac{1 - e^{-\left(\alpha_{i_{j}}(\lambda_{i_{j}}) + 1\right)D_{i_{j}}}}{D_{i_{j}}}\right)\right].\xi_{I,{\star}}\left(ch(Gr^{i_{1}, ..., i_{q}}_{\lambda_{i_{1}}, ..., \lambda_{i_{q}}})\right)
$$
$$
$$

We have  $ch^{Vb}\left(E(-D)\right)e^{D} = ch^{Vb}(E)$, \  therefore:

%%%%%%%%%%%%%%%%%%%%%%%%%%%%%%%%%%%%%%%%%%%%%%%%%%%%%%%%%%%%%%%%%%%%%%%%%%%%%%%%%%%%%%%%%%%%%%%
\begin{theorem}
\label{maintheorem}
\begin{equation}
\label{thhequation}
ch^{Par}(E) = ch^{Vb}(E) \  - 
$$
$$
e^{D}.\sum_{q = 1}^{n} (-1)^{q} \sum_{i_{1}<i_{2}<...<i_{q}}
\sum_{\lambda_{i_{j}} \in \Sigma'_{i_{j}}} \prod_{j = 1}^{q} \left(\frac{1 - e^{-D_{i_{j}}}}{D_{i_{j}}}\right).\xi_{I,{\star}}\left(ch(Gr^{i_{1}, ...,i_{q}}_{\lambda_{i_{1}}, ...,\lambda_{i_{q}}})\right) +
$$
$$
e^{D}.\sum_{q = 1}^{n} (-1)^{q} \sum_{i_{1}<i_{2}<...<i_{q}}
\sum_{\lambda_{i_{j}} \in \Sigma'_{i_{j}}} \prod_{j = 1}^{q} \left[\left(\frac{1 - e^{-\left(\alpha_{i_{j}}(\lambda_{i_{j}}) + 1\right)D_{i_{j}}}}{D_{i_{j}}}\right)\right].\xi_{I,{\star}}\left(ch(Gr^{i_{1}, ..., i_{q}}_{\lambda_{i_{1}}, ..., \lambda_{i_{q}}})\right).
\end{equation}
\end{theorem}
%%%%%%%%%%%%%%%%%%%%%%%%%%%%%%%%%%%%%%%%%%%%%%%%%%%%%%%%%%%%%%%%%%%%%%%%%%%%%%%%%%%%%
%%%%%%%%%%%%%%%%%%%%%%%%%%%%%%%%%%%%%%%%%%%%%%%%%%%%%%%%%%%%%%%%%%%%%%%%%%%%%%%%%%%
$$
$$
\section{Computation of parabolic Chern characters of a locally abelian  parabolic bundle $E$ in codimension one and two
 $ch^{Par}_{1}(E), \ ch^{Par}_{2}(E)$}
$$
$$
\begin{proposition}
For an n-dimensional, non singular variety $Y$, set 
$$
A^{p}Y = A_{n - p} Y,
$$
where $p$ denotes the codimension,  and $n -p$ the dimension.
With this indexing by codimension, the product $x \otimes y \rightarrow  x.y$, reads
$$
A^{p}Y \otimes A^{q}Y \rightarrow A^{p + q}Y,
$$
i.e, the degrees add.
Let $1\in A^{0}Y$ denote the class corresponding to $[Y]$ in $A_{n}Y$, and set
$A^{\ast}Y = \bigoplus A^{p} Y$.
\end{proposition}

%%%%%%%%%%%%%%%%%%%%%%%%%%%%%%%%%%%%%%%%%%%%%%%%%%%%%%%%%%%%%%%%%%%%%%
Return to the equation \eqref{thhequation}
$$
$$
$$
ch^{Par}(E) = ch^{Vb}(E) \  - 
$$
$$
e^{D}.\sum_{q = 1}^{n} (-1)^{q} \sum_{i_{1}<i_{2}<...<i_{q}}
\sum_{\lambda_{i_{j}} \in \Sigma'_{i_{j}}} \prod_{j = 1}^{q} \left(\frac{1 - e^{-D_{i_{j}}}}{D_{i_{j}}}\right).\xi_{I,{\star}}\left(ch(Gr^{i_{1}, ...,i_{q}}_{\lambda_{i_{1}}, ...,\lambda_{i_{q}}})\right) +
$$
$$
e^{D}.\sum_{q = 1}^{n} (-1)^{q} \sum_{i_{1}<i_{2}<...<i_{q}}
\sum_{\lambda_{i_{j}} \in \Sigma'_{i_{j}}} \prod_{j = 1}^{q} \left[\left(\frac{1 - e^{-\left(\alpha_{i_{j}}(\lambda_{i_{j}}) + 1\right)D_{i_{j}}}}{D_{i_{j}}}\right)\right].\xi_{I,{\star}}\left(ch(Gr^{i_{1}, ..., i_{q}}_{\lambda_{i_{1}}, ..., \lambda_{i_{q}}})\right)
$$
$$
$$
take
$$
$$

$\frac{1 - e^{-\left(\alpha_{i_{j}}(\lambda_{i_{j}}) + 1\right)D_{i_{j}}}}{D_{i_{j}}}= 
\frac{1 - \left(1 - (\alpha_{i_{j}}(\lambda_{i_{j}}) + 1)D_{i_{j}} + \frac{(\alpha_{i_{j}}(\lambda_{i_{j}}) + 1)^{2}}{2}D^{2}_{i_{j}} - 
\frac{(\alpha_{i_{j}}(\lambda_{i_{j}}) + 1)^{3}}{6}D^{3}_{i_{j}} + ... \right)}{D_{i_{j}}}$
$$
$$
$$
= \left(\alpha_{i_{j}}(\lambda_{i_{j}}) + 1\right) - \frac{\left(\alpha_{i_{j}}(\lambda_{i_{j}}) + 1\right)^{2}}{2}D_{i_{j}} + 
\frac{\left(\alpha_{i_{j}}(\lambda_{i_{j}}) + 1\right)^{3}}{6}D^{2}_{i_{j}} -...
$$
$$
$$
Let $Gr^{I}_{\lambda}$ be a vector bundles over $D_{I}$  for $ I = (i_{1}, ..., i_{q})$ and $\lambda = (\lambda_{i_{1}}, ..., \lambda_{i_{q}})$ with rank $r$, we have 
$$
$$
\begin{eqnarray*}
\xi_{I,{\star}}\left(ch(Gr^{I}_{\lambda})\right) & = & \xi_{I,{\star}}\left( ch^{D_{I}}_{0}(Gr_{\lambda}^{I}) \right. 
 + ch^{D_{I}}_{1}(Gr^{I}_{\lambda}) \\
 && + ch^{D_{I}}_{2}(Gr^{I}_{\lambda}) 
  + ch^{D_{I}}_{3}(Gr^{I}_{\lambda}) 
  + \left.  ...\right)/X \\
 & =  & \xi_{I,{\star}}\left(ch^{D_{I}}_{0}(Gr_{\lambda}^{I})\right)/X 
  + \xi_{I,{\star}}\left(ch^{D_{I}}_{1}(Gr^{I}_{\lambda})\right)/X \\
&& + \xi_{I,{\star}}\left(ch^{D_{I}}_{2}(Gr^{I}_{\lambda})\right)/X   
 + \xi_{I,{\star}}\left(ch^{D_{I}}_{3}(Gr^{I}_{\lambda})\right)/X 
 + ... .
\end{eqnarray*}
so 
$$
ch^{D_{I}}_{0}(Gr^{I}_{\lambda})  =  (rank(Gr^{I}_{\lambda}))/D_{I} 
 = 
\sum_{p\in Irr(D_{I})}rank_{p}(Gr^{I}_{\lambda}).[D_{p}]/D_{I} 
$$
$$ \mbox{ in  }  A^{0}(D_{I}) = A_{dimD_{I}}(D_{I}) 
 = \bigoplus A^{0}(D_{p}) = \bigoplus \mathbb{Q}.[D_{p}]
$$
where $D_{I} = \bigcup_{p \in Irr(D_{I})}D_{p}$ 
and $Irr(D_{I})$ denotes the set of irreducible components of 
$D_{i_{1}} \cap D_{i_{2}} \cap ... \cap D_{i_{q}}$ then
$$
$$ 
$\xi_{I,{\star}}\left(\sum_{p\in Irr(D_{I})}rank_{p}(Gr^{I}_{\lambda}).[D_{p}]/D_{I}\right) = \left(\sum_{p\in Irr(D_{I})}rank_{p}(Gr^{I}_{\lambda}).[D_{p}]/X\right) \in A^{q}(X)$ \ \ \ \  
$$
$$
which is of codimension $q$,
$$
$$
$ch^{D_{I}}_{1}(Gr^{I}_{\lambda})= c^{D_{I}}_{1}(Gr^{I}_{\lambda}) \in A^{1}(D_{I}) = A_{dimD_{I}-1}(D_{I})$
then
$$
$$
$\left(\xi_{I,{\star}}\left(c^{D_{I}}_{1}(Gr^{I}_{\lambda})\right)/X \right) \in A^{1 + q}(X)$ which is  of codimension $q + 1$,
$$
$$
$ch^{D_{I}}_{2}(Gr^{I}_{\lambda}) = \frac{1}{2}\left[(c^{D_{I}}_{1})^{2}(Gr^{I}_{\lambda})- 2c_{2}^{D_{I}}(Gr^{I}_{\lambda})\right] \in A^{2}(D_{I}) = A_{dimD_{I} - 2}(D_{I})$ then
$$
$$
$ \xi_{I,{\star}}\left(\frac{1}{2}[(c^{D_{I}}_{1})^{2}(Gr^{I}_{\lambda})- 2c_{2}^{D_{I}}(Gr^{I}_{\lambda})]\right)/X \in A^{2 + q}$
which is of codimension $q + 2$, 
$$
$$
$ch_{3}^{D_{I}}(Gr^{I}_{\lambda}) =  \frac{1}{6}\left[(c^{D_{I}}_{1})^{3}(Gr^{I}_{\lambda})  - 3c^{D_{I}}_{1}(Gr^{I}_{\lambda})c^{D_{I}}_{2}(Gr^{I}_{\lambda})] + 3c^{D_{I}}_{3}(Gr^{I}_{\lambda})\right] \in A^{3}(D_{I})$  then
$$
$$
$ \xi_{I,{\star}} \left(\frac{1}{6}\left[(c^{D_{I}}_{1})^{3}(Gr^{I}_{\lambda}) - 3(c^{D_{I}}_{1}(Gr^{I}_{\lambda})c^{D_{I}}_{2}(Gr^{I}_{\lambda}) + 3c^{D_{I}}_{3}(Gr^{I}_{\lambda}) \right]\right)/X \in A^{3 + q}$  
$$
$$
which is of codimension  $q + 3$. Therefore
$$
$$
\begin{equation}
ch^{Par}(E) = ch^{Vb}(E) + e^{D}.\sum_{q = 1}^{n} (-1)^{q} \sum_{i_{1}<i_{2}<...<i_{q}}
\sum_{\lambda_{i_{j}} \in \Sigma'_{i_{j}}}
$$
$$
\prod_{j = 1}^{q}\left[\left(\alpha_{i_{j}}(\lambda_{i_{j}}) + 1\right) - \frac{\left(\alpha_{i_{j}}(\lambda_{i_{j}}) + 1\right)^{2}}{2}D_{i_{j}} + 
\frac{\left(\alpha_{i_{j}}(\lambda_{i_{j}}) + 1\right)^{3}}{6}D^{2}_{i_{j}} -...\right].
$$
$$
(\sum_{p\in Irr(D_{I})}rank_{p}(Gr^{I}_{\lambda}).[D_{p}]/X + \xi_{I,{\star}}(c^{D_{I}}_{1}(Gr^{I}_{\lambda}))/X
+ \xi_{I,{\star}}(\frac{1}{2}[(c^{D_{I}}_{1})^{2}(Gr^{I}_{\lambda})- 2c_{2}^{D_{I}}(Gr^{I}_{\lambda})])/X
$$
$$
+ \xi_{I,{\star}} (\frac{1}{6}[(c^{D_{I}}_{1})^{3}(Gr^{I}_{\lambda}) - 3(c^{D_{I}}_{1}(Gr^{I}_{\lambda})c^{D_{I}}_{2}(Gr^{I}_{\lambda}) +
3c^{D_{I}}_{3}(Gr^{I}_{\lambda})])) \ 
$$
$$
- e^{D}.\sum_{q = 1}^{n} (-1)^{q} \sum_{i_{1}<i_{2}<...<i_{q}}
\sum_{\lambda_{i_{j}} \in \Sigma'_{i_{j}}} \prod_{j = 1}^{q}\left( 1 - \frac{D_{i_{j}}}{2} + \frac{D^{2}_{i_{j}}}{6} -...\right) \cdot
\left[ \sum_{p\in Irr(D_{I})}r_{p}(Gr^{I}_{\lambda}).[D_{p}]/X + \right.
$$
$$
\xi_{I,{\star}}(c^{D_{I}}_{1}(Gr^{I}_{\lambda}))/X + \xi_{I,{\star}}(\frac{1}{2}[(c^{D_{I}}_{1})^{2}(Gr^{I}_{\lambda})- 2c_{2}^{D_{I}}(Gr^{I}_{\lambda})])/X + 
$$
$$
\left. \xi_{I,{\star}} (\frac{1}{6}[(c^{D_{I}}_{1})^{3}(Gr^{I}_{\lambda}) - 3(c^{D_{I}}_{1}(Gr^{I}_{\lambda})c^{D_{I}}_{2}(Gr^{I}_{\lambda}) 
+ 3c^{D_{I}}_{3}(Gr^{I}_{\lambda})])\right]
\end{equation}
$$
$$
we have 
$$
$$
$ch^{Vb}(E) = rank(E).[X] + ch^{Vb}_{1}(E) + ch_{2}^{Vb}(E) + ch^{Vb}_{3}(E) + ...$
$$
$$
where  \ $rank(E).[X] \in A^{0}, \ ch_{1}^{Vb}(E) \in A^{1},\   ch_{2}^{Vb}(E) \in A^{2},\   ch_{3}^{Vb}(E) \in A^{3}$,
$$
$$
$ and \ e^{D} = 1 + D + \frac{D^{2}}{2} + \frac{D^{3}}{6} + ...$  \  \ where  $1 \in A^{0},\  D \in A^{1}, \  \frac{D^{2}}{2} \in A^{2},\   \frac{D^{3}}{6} \in A^{3}$
$$
$$
for $k = 0,1, ..., n$, \  $ch^{Par}_{k}(E) = ch_{0}^{Par}(E)  + ch_{1}^{Par}(E) + ch_{2}^{Par}(E) + ch_{3}^{Par}(E) + ...$ 
$$
$$
then
$$
$$
$ch_{0}^{Par}(E)  + ch_{1}^{Par}(E) + ch_{2}^{Par}(E)  = rank(E).[X]/A^{0} + ch^{Vb}_{1}(E)/A^{1}$ 
$$
$$
$ + \ ch^{Vb}_{2}(E)/A^{2} - \sum_{i_{1} \in \mathcal{S}}\sum_{\lambda_{i_{1}} \in \sum'_{i_{1}}}(\alpha_{i_{1}}(\lambda_{i_{1}}) + 1 ).rank(Gr^{i_{1}}_{\lambda_{i_{1}}}).[D_{i_{1}}]/A^{1}$
$$
$$
$ - \  \sum_{i_{1} \in \mathcal{S}}\sum_{\lambda_{i_{1}} \in \sum'_{i_{1}}}(\alpha_{i_{1}}(\lambda_{i_{1}}) + 1).(\xi_{i_{1}})_{\star}\left(c_{1}^{D_{i_{1}}}(Gr^{i_{1}}_{\lambda_{i_{1}}})\right)/A^{2}$
$$
$$
$ - \ \sum_{i_{1} \in \mathcal{S}}\sum_{\lambda_{i_{1}} \in \sum'_{i_{1}}} (\alpha_{i_{1}}(\lambda_{i_{1}}) + 1 ).rank(Gr^{i_{1}}_{\lambda_{i_{1}}}).[D_{i_{1}}].[D]/A^{2}$ 
$$
$$
$ + \sum_{i_{1} <i_{2}}
 \sum_{\aatop{\lambda_{i_{1}} }{ \lambda_{i_{2}}}}  \sum_{p \in Irr(D_{i_{1}} \cap D_{i_{2}})}(\alpha_{i_{1}}(\lambda_{i_{1}}) + 1 )(\alpha_{i_{2}}(\lambda_{i_{2}}) + 1 ).rank_{p}(Gr^{i_{1},i_{2}}_{\lambda_{i_{1}},\lambda_{i_{2}}}).[D_{p}]/A^{2}$
 $$
 $$
$ + \frac{1}{2}\sum_{i_{1} \in \mathcal{S}}\sum_{\lambda_{i_{1}} \in \sum'_{i_{1}}}(\alpha_{i_{1}}(\lambda_{i_{1}}) + 1)^{2}.rank(Gr^{i_{1}}_{\lambda_{i_{1}}}).[D_{i_{1}}]^{2}/A^{2}$
$$
$$
%%%%%%%%%%%%%%%%%%%%%%%%%%%%%%%%%%%%%%%%%%%%%%%%%%%%%%%%%%%%%%%%%%%%%%%%%
$ + \sum_{i_{1} \in \mathcal{S}}\sum_{\lambda_{i_{1}} \in \sum'_{i_{1}}} rank(Gr^{i_{1}}_{\lambda_{i_{1}}}).[D_{i_{1}}]/A^{1}$
$$
$$
$ + \sum_{i_{1} \in \mathcal{S}}\sum_{\lambda_{i_{1}} \in \sum'_{i_{1}}} (\xi_{i_{1}})_{\star}\left(c_{1}^{D_{i_{1}}}(Gr^{i_{1}}_{\lambda_{i_{1}}})\right)/A^{2}$
$$
$$
$ + \sum_{i_{1} \in \mathcal{S}}\sum_{\lambda_{i_{1}} \in \sum'_{i_{1}}} rank(Gr^{i_{1}}_{\lambda_{i_{1}}}).[D_{i_{1}}].[D]/A^{2}$ 
$$
$$
$ - \sum_{i_{1} <i_{2}}
 \sum_{\aatop{\lambda_{i_{1}} }{ \lambda_{i_{2}}}}  \sum_{p \in Irr(D_{i_{1}} \cap D_{i_{2}})} rank_{p}(Gr^{i_{1},i_{2}}_{\lambda_{i_{1}},\lambda_{i_{2}}}).[D_{p}]/A^{2}$
 $$
 $$
$ - \  \frac{1}{2}\sum_{i_{1} \in \mathcal{S}}\sum_{\lambda_{i_{1}} \in \sum'_{i_{1}}} rank(Gr^{i_{1}}_{\lambda_{i_{1}}}).[D_{i_{1}}]^{2}/A^{2}$
$$
$$
%%%%%%%%%%%%%%%%%%%%%%%%%%%%%%%%%%%%%%%%%%%%%%%%%%%%%%%%%%%%%%%%%%%%%%%%%%%%%%%%%%%%%%%%%%%%%%%%
\begin{lemma}
In the Grothendieck group, for all $\ i_{1}, i_{2} \in \mathcal{S}$, for all 
$\lambda_{i_{1}} \in \Sigma'_{i_{1}}$, for all $\ p \in Irr(D_{i_{1}} \cap D_{i_{2}})$, we have:
$$
rank(Gr^{i_{1}}_{\lambda_{i_{1}}}) = \sum_{\lambda_{i_{2}} \in \Sigma'_{i_{2}}}rank_{p}(Gr^{i_{1}, i_{2}}_{\lambda_{i_{1}}, \lambda_{i_{2}}}).
$$
\end{lemma}
%%%%%%%%%%%%%%%%%%%%%%%%%%%%%%%%%%%%%%%%%%%%%%%%%%%%%%%%%%%%%%%%%%%%%%%%%%%%%%%%%%%%%%%%%%%%%%%%%
$$
$$
Now we have  $D = \sum_{i} D_{i}$ \  and $[D_{_{i_{1}}}].[D_{i_{2}}] = \sum_{p \in Irr(D_{i_{1}} \cap D_{i_{2}})}[D_{p}]$ \ then 
$$
$$
$[D].[D_{i_{1}}] = \sum_{i_{1} \in \mathcal{S}}[D_{i_{1}}]^{2} + \sum_{i_{1} \neq i_{2}}[D_{i_{1}}].[D_{i_{2}}] =  
\sum_{i_{1} \in \mathcal{S}}[D_{i_{1}}]^{2} + \sum_{i_{1} \neq i_{2}}\sum_{p \in Irr(D_{i_{1}} \cap D_{i_{2}})}[D_{p}].$
$$
$$
Then
$$
$$
$ch_{0}^{Par}(E)  + ch_{1}^{Par}(E) + ch_{2}^{Par}(E)  = rank(E).[X]/A^{0} + ch^{Vb}_{1}(E)/A^{1}$ 
$$
$$
$ + \  ch^{Vb}_{2}(E)/A^{2}  - \sum_{i_{1} \in \mathcal{S}}
\sum_{\lambda_{i_{1}} \in \sum'_{i_{1}}} \alpha_{i_{1}}(\lambda_{i_{1}}).rank(Gr^{i_{1}}_{\lambda_{i_{1}}}).[D_{i_{1}}]/A^{1}$
$$
$$
$ - \sum_{i_{1} \in \mathcal{S}}\sum_{\lambda_{i_{1}} \in \sum'_{i_{1}}}\alpha_{i_{1}}(\lambda_{i_{1}}).(\xi_{i_{1}})_{\star}\left(c_{1}^{D_{i_{1}}}(Gr^{i_{1}}_{\lambda_{i_{1}}})\right)/A^{2}$
$$
$$
$ - \sum_{i_{1} \in \mathcal{S}}\sum_{\lambda_{i_{1}} \in \sum'_{i_{1}}} \alpha_{i_{1}}(\lambda_{i_{1}}) .rank(Gr^{i_{1}}_{\lambda_{i_{1}}}).[D_{i_{1}}]^{2}/A^{2}$  
$$
$$
$ - \sum_{i_{1} \in \mathcal{S}} \sum_{i_{1} \neq i_{2}}\sum_{\lambda_{i_{1}} \in \sum'_{i_{1}}} \sum_{\lambda_{i_{2}} \in \Sigma'_{i_{2}}} \sum_{p \in Irr(D_{i_{1}} \cap D_{i_{2}})}\alpha_{i_{1}}(\lambda_{i_{1}}).rank_{p}(Gr^{i_{1},i_{2}}_{\lambda_{i_{1}},\lambda_{i_{2}}}).[D_{p}]/A^{2}$
$$
$$
$ + \sum_{i_{1} <i_{2}}
 \sum_{\aatop{\lambda_{i_{1}} }{ \lambda_{i_{2}}}}  \sum_{p \in Irr(D_{i_{1}} \cap D_{i_{2}})} \alpha_{i_{1}}(\lambda_{i_{1}}).rank_{p}(Gr^{i_{1}i_{2}}_{\lambda_{i_{1}},\lambda_{i_{2}}}).[D_{p}]/A^{2}$
$$
$$
$ + \sum_{i_{1} <i_{2}}
\sum_{\aatop{\lambda_{i_{1}} }{ \lambda_{i_{2}}}}  \sum_{p \in Irr(D_{i_{1}} \cap D_{i_{2}})} \alpha_{i_{2}}(\lambda_{i_{2}}).rank_{p}(Gr^{i_{1}i_{2}}_{\lambda_{i_{1}},\lambda_{i_{2}}}).[D_{p}]/A^{2}$
$$
$$
$ + \sum_{i_{1} <i_{2}}
\sum_{\aatop{\lambda_{i_{1}} }{ \lambda_{i_{2}}}}  \sum_{p \in Irr(D_{i_{1}} \cap D_{i_{2}})} \alpha_{i_{1}}(\lambda_{i_{1}}).\alpha_{i_{2}}(\lambda_{i_{2}}).rank_{p}(Gr^{i_{1}i_{2}}_{\lambda_{i_{1}},\lambda_{i_{2}}}).[D_{p}]/A^{2}$
$$
$$
$ +  \ \frac{1}{2}\sum_{i_{1} \in \mathcal{S}}\sum_{\lambda_{i_{1}} \in \sum'_{i_{1}}} \alpha^{2}_{i_{1}}(\lambda_{i_{1}}).rank(Gr^{i_{1}}_{\lambda_{i_{1}}}).[D_{i_{1}}]^{2}/A^{2}$
$$
$$
$ + \sum_{i_{1} \in \mathcal{S}}\sum_{\lambda_{i_{1}} \in \sum'_{i_{1}}} \alpha_{i_{1}}(\lambda_{i_{1}}).rank(Gr^{i_{1}}_{\lambda_{i_{1}}}).[D_{i_{1}}]^{2}/A^{2}$
$$
$$
%%%%%%%%%%%%%%%%%%%%%%%%%%%%%%%%%%%%%%%%%%%%%%%%%%%%%%%%%%%%%%%%%%%%%%%%%%%%%%%%%%%%%%%%%%%%%%%%%%%%%%%%%%
%%%%%%%%%%%%%%%%%%%%%%%%%%%%%%%%%%%%%%%%%%%%%%%%%%%%%%%%%%%%%%%%%%%%%%%%%%%%%%%%%%%%%%%%%%%%%%%%%%%%%%%%%%
\subsection{The characteristic numbers for parabolic bundle in codimension one and  two.}
-For any parabolic bundle $E$  in codimension one, and two,  the parabolic first, second  Chern characters 
$ch_{1}^{Par}(E), and \  ch_{2}^{Par}(E),\ $  are obtained as follows:
$$
$$
$\centerdot  \ ch_{0}^{Par}(E) := rank(E).[X]$
$$
$$
$\centerdot \  ch^{Par}_{1}(E) := ch_{1}^{Vb}(E) \ - \ \sum_{i_{1} \in \mathcal{S}}\sum_{\lambda_{i_{1}} \in 
\Sigma '_{i_{1}}}\alpha_{i_{1}}(\lambda_{i_{1}}).rank(Gr^{i_{1}}_{\lambda_{i_{1}}}).[D_{i_{1}}]$
$$
$$
$\centerdot \ ch^{Par}_{2}(E) := \  ch^{Vb}_{2}(E) \ -\  \sum_{i_{1} \in \mathcal{S}}\sum_{\lambda_{i_{1}} \in \sum'_{i_{1}}} \alpha_{i_{1}}(\lambda_{i_{1}}).(\xi_{i_{1}})_{\star}\left(c_{1}^{D_{i_{1}}}(Gr^{i_{1}}_{\lambda_{i_{1}}})\right) $
$$
$$
$\hspace{2.2cm} +  \ \frac{1}{2} \ \sum_{i_{1} \in \mathcal{S}}\sum_{\lambda_{i_{1}} \in \sum'_{i_{1}}} \alpha^{2}_{i_{1}}(\lambda_{i_{1}}).rank(Gr^{i_{1}}_{\lambda_{i_{1}}}).[D_{i_{1}}]^{2}$
$$
$$
$\hspace{2.2cm} +  \sum_{i_{1} <i_{2}}
\sum_{\aatop{\lambda_{i_{1}} }{ \lambda_{i_{2}}}}  \sum_{p \in Irr(D_{i_{1}} \cap D_{i_{2}})} \alpha_{i_{1}}(\lambda_{i_{1}}).\alpha_{i_{2}}(\lambda_{i_{2}}).rank_{p}(Gr^{i_{1}, i_{2}}_{\lambda_{i_{1}}, \lambda_{i_{2}}}).[D_{p}].$
$$
$$
In order to compare with Mochizuki's formula, note that 
$$
$$
$\sum_{i_{1} <i_{2}}
\sum_{\aatop{\lambda_{i_{1}} }{ \lambda_{i_{2}}}}  \sum_{p \in Irr(D_{i_{1}} \cap D_{i_{2}})} \alpha_{i_{1}}(\lambda_{i_{1}}).\alpha_{i_{2}}(\lambda_{i_{2}}).rank_{p}(Gr^{i_{1}, i_{2}}_{\lambda_{i_{1}}, \lambda_{i_{2}}}).[D_{p}] =  $ 
$$
$$
$ \frac{1}{2}\sum_{i_{1} \neq i_{2}}
\sum_{\aatop{\lambda_{i_{1}} }{ \lambda_{i_{2}}}}  \sum_{p \in Irr(D_{i_{1}} \cap D_{i_{2}})} \alpha_{i_{1}}(\lambda_{i_{1}}).\alpha_{i_{2}}(\lambda_{i_{2}}).rank_{p}(Gr^{i_{1}, i_{2}}_{\lambda_{i_{1}}, \lambda_{i_{2}}}).[D_{p}]$
$$
$$
therefore our formula may be written:
$$
$$
$\centerdot \ ch^{Par}_{2}(E) := \  ch^{Vb}_{2}(E) \ -\  \sum_{i_{1} \in \mathcal{S}}\sum_{\lambda_{i_{1}} \in \sum'_{i_{1}}} \alpha_{i_{1}}(\lambda_{i_{1}}).(\xi_{i_{1}})_{\star}\left(c_{1}^{D_{i_{1}}}(Gr^{i_{1}}_{\lambda_{i_{1}}})\right) $
$$
$$
$\hspace{2.2cm} +  \ \frac{1}{2} \ \sum_{i_{1} \in \mathcal{S}}\sum_{\lambda_{i_{1}} \in \sum'_{i_{1}}} \alpha^{2}_{i_{1}}(\lambda_{i_{1}}).rank(Gr^{i_{1}}_{\lambda_{i_{1}}}).[D_{i_{1}}]^{2}$
$$
$$
$ \hspace{2.2cm} +  \ \frac{1}{2}\sum_{i_{1} \neq i_{2}}
\sum_{\aatop{\lambda_{i_{1}} }{ \lambda_{i_{2}}}}  \sum_{p \in Irr(D_{i_{1}} \cap D_{i_{2}})} \alpha_{i_{1}}(\lambda_{i_{1}}).\alpha_{i_{2}}(\lambda_{i_{2}}).rank_{p}(Gr^{i_{1}, i_{2}}_{\lambda_{i_{1}}, \lambda_{i_{2}}}).[D_{p}].$

This coincides exactly with the formula given by Mochizuki in \cite[\S 3.1.5, p. 30]{Mochizuki}. It is also the
same as the definition given by Panov \cite{Panov}. \label{panovpage}
Note that in Panov's general definition the sum for the last term is written as ${\scriptstyle \sum _{i,j}}$ without the factor of $1/2$ 
but later he uses it as a sum over $i<j$, so our formula and Mochizuki's also coincide with Panov's formula
in the way he uses it. 

Mochizuki's formula was for the Chern character in cohomology, which he defined as the 
integral of the Chern form of the curvature of an adapted metric. Our calculation verifies that this gives
the same answer as the method using Deligne-Mumford stacks of \cite{Biswas1} \cite{Borne} 
\cite{IyerSimpson1} for rational weights. 
Our formula is valid for the Chern character in the rational or real Chow ring.

$$
$$

Here we explain some of the notation:
$$
$$
$\centerdot \  ch_{1}^{Vb}(E), ch^{Vb}_{2}(E)$ denotes the first, second, Chern character of vector bundles E.
$$
$$
$\centerdot \  Irr(D_{I})$ denotes the set of the irreducible components of $D_I:=D_{i_{1}} \cap D_{i_{2}}  \cap ... \cap D_{i_{q}}$.
$$
$$
$\centerdot \  \xi_{I}$ denotes the closed immersion $D_{I} \longrightarrow X$, and $\xi_{I,{\star}} : A^k(D_{I}) \longrightarrow A^{k+q}(X)$ denotes  the associated Gysin map.
$$
$$
$\centerdot$ \   Let $p$ be an element of $Irr(D_{i} \cap D_{j}).$  Then $rank_{p}(Gr^{I}_{\lambda})$ denotes the rank of $Gr^{I}_{\lambda}$ as 
an $\mathcal{O}_{p}$-module.
$$
$$
$\centerdot  \  [D_{i_j}] \in A^1(X)\otimes \mathbb{Q}$,  and $[D_{p}] \in A^2(X)\otimes \mathbb{Q}$ denote the cycle classes given by $D_{i_j}$ 
and $D_{p}$ respectively.
%%%%%%%%%%%%%%%%%%%%%%%%%%%%%%%%%%%%%%%%%%%%%%%%%%%%%%%%%%%%%%%%%%%%%%%%%%%%%%%%%%%%%%%%%%%%%%%%%%%%%%%%%%%%%%%%%
%%%%%%%%%%%%%%%%%%%%%%%%%%%%%%%%%%%%%%%%%%%%%%%%%%%%%%%%%%%%%%%%%%%%%%%%%%%%%%%%%%%%%%%%%%%%%%%%%%%%%%%%%%%%%%%%%
$$
$$
\section{Parabolic Chern character of a locally abelian parabolic bundle $E$ in codimension 3, 
$Ch_{3}^{Par}(E)$}
$$
$$
By the same method of computation as above, we get the following formula, which has not been considered elsewhere in the literature. 
$$
$$
$ ch_{3}^{Par}(E) =  ch^{Vb}_{3}(E) - \frac{1}{2}\sum_{i_{1} \in \mathcal{S}}\sum_{\lambda_{i_{1}} \in \sum'_{i_{1}}} (\alpha_{i_{1}}(\lambda_{i_{1}}) + 1 ).rank(Gr^{i_{1}}_{\lambda_{i_{1}}}).[D_{i_{1}}].[D]^{2}$
$$
$$
$ + \frac{1}{2}\sum_{i_{1} \in \mathcal{S}}\sum_{\lambda_{i_{1}} \in \sum'_{i_{1}}} (\alpha_{i_{1}}(\lambda_{i_{1}}) + 1 )^{2}.(\xi_{i_{1}})_{\star}\left(c_{1}^{D_{i_{1}}}(Gr^{i_{1}}_{\lambda_{i_{1}}})\right).[D_{i_{1}}]$
$$
$$
$ - \frac{1}{6} \sum_{i_{1} \in \mathcal{S}}\sum_{\lambda_{i_{1}} \in \sum'_{i_{1}}} (\alpha_{i_{1}}(\lambda_{i_{1}}) + 1)^{3}.rank(Gr^{i_{1}}_{\lambda_{i_{1}}}).[D_{i_{1}}]^{3}$
$$
$$
$ - \sum_{i_{1} \in \mathcal{S}}\sum_{\lambda_{i_{1}} \in \sum'_{i_{1}}} (\alpha_{i_{1}}(\lambda_{i_{1}}) + 1 ).(\xi_{i_{1}})_{\star}\left(c_{1}^{D_{i_{1}}}(Gr^{i_{1}}_{\lambda_{i_{1}}})\right).[D]$
$$
$$
$ + \frac{1}{2}\sum_{i_{1} \in \mathcal{S}}
\sum_{\lambda_{i_{1}} \in \sum'_{i_{1}}} (\alpha_{i_{1}}(\lambda_{i_{1}}) + 1 )^{2}.rank(Gr^{i_{1}}_{\lambda_{i_{1}}}).[D_{i_{1}}]^{2}.[D]$
$$
$$
$ - \frac{1}{2}\sum_{i_{1} \in \mathcal{S}}\sum_{\lambda_{i_{1}} \in \sum'_{i_{1}}}(\alpha_{i_{1}}(\lambda_{i_{1}}) + 1 ).(\xi_{i_{1}})_{\star}\left((c_{1}^{D_{i_{1}}})^{2}(Gr^{i_{1}}_{\lambda_{i_{1}}}) - 2c_{2}^{D_{i_{1}}}(Gr_{\lambda_{i_{1}}}^{i_{1}})\right)$
$$
$$
$ - \frac{1}{2} \sum_{i_{1} <i_{2}} \sum_{\aatop{\lambda_{i_{1}} }{ \lambda_{i_{2}}}}  \sum_{p \in Irr(D_{i_{1}} \cap D_{i_{2}})}(\alpha_{i_{1}}(\lambda_{i_{1}}) + 1 ).(\alpha_{i_{2}}(\lambda_{i_{2}}) + 1 )^{2}.rank_{p}(Gr^{i_{1},i_{2}}_{\lambda_{i_{1}},\lambda_{i_{2}}}).[D_{i_{2}}].[D_{p}]$ 
$$
$$
$ - \frac{1}{2} \sum_{i_{1} <i_{2}} \sum_{\aatop{\lambda_{i_{1}} }{ \lambda_{i_{2}}}}  \sum_{p \in Irr(D_{i_{1}} \cap D_{i_{2}})}(\alpha_{i_{1}}(\lambda_{i_{1}}) + 1 )^{2}.(\alpha_{i_{2}}(\lambda_{i_{2}}) + 1 ).rank_{p}(Gr^{i_{1},i_{2}}_{\lambda_{i_{1}},\lambda_{i_{2}}}).[D_{i_{1}}].[D_{p}]$ 
$$
$$
$ + \sum_{i_{1} <i_{2}} \sum_{\aatop{\lambda_{i_{1}} }{ \lambda_{i_{2}}}}  \sum_{p \in Irr(D_{i_{1}} \cap D_{i_{2}})}(\alpha_{i_{1}}(\lambda_{i_{1}}) + 1 ).(\alpha_{i_{2}}(\lambda_{i_{2}}) + 1 ).rank_{p}(Gr^{i_{1},i_{2}}_{\lambda_{i_{1}},\lambda_{i_{2}}}).[D].[D_{p}]$ 
$$
$$
$+ \sum_{i_{1} <i_{2}} \sum_{\aatop{\lambda_{i_{1}} }{ \lambda_{i_{2}}}}(\alpha_{i_{1}}(\lambda_{i_{1}}) + 1 ).(\alpha_{i_{2}}(\lambda_{i_{2}}) + 1 ).(\xi_{i_{1},i_{2}})_{\star}\left(c_{1}^{D_{i_{1}} \cap D_{i_{2}}}(Gr^{i_{1},i_{2}}_{\lambda_{i_{1}},\lambda_{i_{2}}})\right)$
$$
$$
$ - \sum_{i_{1} <i_{2} < i_{3}} \sum_{\aatop{\lambda_{i_{1}} }{  \lambda_{i_{3}}}}  \sum_{p \in Irr(D_{i_{1}} \cap D_{i_{2}} \cap D_{i_{3}})}(\alpha_{i_{1}}(\lambda_{i_{1}}) + 1 ).(\alpha_{i_{2}}(\lambda_{i_{2}}) + 1 ).(\alpha_{i_{3}}(\lambda_{i_{3}}) + 1 ).rank_{p}(Gr^{i_{1},i_{2},i_{3}}_{\lambda_{i_{1}},\lambda_{i_{2}},\lambda_{i_{3}}}).[D_{p}]$
$$
$$
$ + \frac{1}{2}\sum_{i_{1} \in \mathcal{S}}\sum_{\lambda_{i_{1}} \in \sum'_{i_{1}}} rank(Gr^{i_{1}}_{\lambda_{i_{1}}}).[D_{i_{1}}].[D]^{2}
 - \frac{1}{2}\sum_{i_{1} \in \mathcal{S}}\sum_{\lambda_{i_{1}} \in \sum'_{i_{1}}} (\xi_{i_{1}})_{\star}\left(c_{1}^{D_{i_{1}}}(Gr^{i_{1}}_{\lambda_{i_{1}}})\right).[D_{i_{1}}]$
$$
$$
$ \ + \frac{1}{6} \sum_{i_{1} \in \mathcal{S}}\sum_{\lambda_{i_{1}} \in \sum'_{i_{1}}} rank(Gr^{i_{1}}_{\lambda_{i_{1}}}).[D_{i_{1}}]^{3}
+ \sum_{i_{1} \in \mathcal{S}}\sum_{\lambda_{i_{1}} \in \sum'_{i_{1}}} (\xi_{i_{1}})_{\star}\left(c_{1}^{D_{i_{1}}}(Gr^{i_{1}}_{\lambda_{i_{1}}})\right).[D]$
$$
$$
$ - \frac{1}{2}\sum_{i_{1} \in \mathcal{S}}\sum_{\lambda_{i_{1}} \in \sum'_{i_{1}}} rank(Gr^{i_{1}}_{\lambda_{i_{1}}}).[D_{i_{1}}]^{2}.[D]
\  + \ \frac{1}{2}\sum_{i_{1} \in \mathcal{S}}\sum_{\lambda_{i_{1}} \in \sum'_{i_{1}}} (\xi_{i_{1}})_{\star}\left((c_{1}^{D_{i_{1}}})^{2}(Gr^{i_{1}}_{\lambda_{i_{1}}}) - 2c_{2}^{D_{i_{1}}}(Gr_{\lambda_{i_{1}}}^{i_{1}})\right)$
$$
$$
$ - \sum_{i_{1} <i_{2}} \sum_{\aatop{\lambda_{i_{1}} }{ \lambda_{i_{2}}}}  \sum_{p \in Irr(D_{i_{1}} \cap D_{i_{2}})} rank_{p}(Gr^{i_{1},i_{2}}_{\lambda_{i_{1}},\lambda_{i_{2}}}).[D].[D_{p}] \ - \sum_{i_{1} <i_{2}} \sum_{\aatop{\lambda_{i_{1}} }{ \lambda_{i_{2}}}} (\xi_{i_{1},i_{2}})_{\star}\left(c_{1}^{D_{i_{1}} \cap D_{i_{2}}}(Gr^{i_{1},i_{2}}_{\lambda_{i_{1}},\lambda_{i_{2}}})\right)$
$$
$$
$ + \sum_{i_{1} <i_{2} < i_{3}} \sum_{\aatop{\lambda_{i_{1}} }{  \lambda_{i_{3}}}}  \sum_{p \in Irr(D_{i_{1}} \cap D_{i_{2}} \cap D_{i_{3}})} rank_{p}(Gr^{i_{1},i_{2},i_{3}}_{\lambda_{i_{1}},\lambda_{i_{2}},\lambda_{i_{3}}}).[D_{p}]$
$$
$$
$ + \frac{1}{2}\sum_{i_{1} <i_{2}} \sum_{\aatop{\lambda_{i_{1}} }{ \lambda_{i_{2}}}}  \sum_{p \in Irr(D_{i_{1}} \cap D_{i_{2}})} rank_{p}(Gr^{i_{1},i_{2}}_{\lambda_{i_{1}},\lambda_{i_{2}}}).[D_{i_{2}}].[D_{p}]$ 
$$
$$
$ + \frac{1}{2}\sum_{i_{1} <i_{2}} \sum_{\aatop{\lambda_{i_{1}} }{ \lambda_{i_{2}}}}  \sum_{p \in Irr(D_{i_{1}} \cap D_{i_{2}})} rank_{p}(Gr^{i_{1},i_{2}}_{\lambda_{i_{1}},\lambda_{i_{1}}}).[D_{i_{1}}].[D_{p}]$
$$
$$
\subsection{The characteristic number for parabolic bundle in codimension three.}
-For any parabolic bundle $E$  in codimension 3,  the parabolic third Chern character $\  ch_{3}^{Par}(E),\ $  is obtained as follows:
$$
$$
$ ch_{3}^{Par}(E) =  ch^{Vb}_{3}(E) - \frac{1}{2}\sum_{i_{1} \in \mathcal{S}}\sum_{\lambda_{i_{1}} \in \sum'_{i_{1}}} \alpha_{i_{1}}(\lambda_{i_{1}}).rank(Gr^{i_{1}}_{\lambda_{i_{1}}}).[D_{i_{1}}].[D]^{2}$
$$
$$
$\hspace{.6cm} + \ \frac{1}{2}\sum_{i_{1} \in \mathcal{S}}\sum_{\lambda_{i_{1}} \in \sum'_{i_{1}}} \left[\alpha^{2}_{i_{1}}(\lambda_{i_{1}}) + 2 \alpha_{i_{1}}(\lambda_{i_{1}})\right].(\xi_{i_{1}})_{\star}\left(c_{1}^{D_{i_{1}}}(Gr^{i_{1}}_{\lambda_{i_{1}}})\right).[D_{i_{1}}]$
$$
$$
$ \hspace{.6cm} - \frac{1}{6} \sum_{i_{1} \in \mathcal{S}}\sum_{\lambda_{i_{1}} \in \sum'_{i_{1}}} [\alpha^{3}_{i_{1}}(\lambda_{i_{1}}) + 3 \alpha^{2}_{i_{1}}(\lambda_{i_{1}}) + 3 \alpha_{i_{1}}(\lambda_{i_{1}})].rank(Gr^{i_{1}}_{\lambda_{i_{1}}}).[D_{i_{1}}]^{3}$
$$
$$
$\hspace{.6cm} - \sum_{i_{1} \in \mathcal{S}}\sum_{\lambda_{i_{1}} \in \sum'_{i_{1}}} \alpha_{i_{1}}(\lambda_{i_{1}}).(\xi_{i_{1}})_{\star}\left(c_{1}^{D_{i_{1}}}(Gr^{i_{1}}_{\lambda_{i_{1}}})\right).[D]$
$$
$$
$\hspace{.6cm} + \frac{1}{2}\sum_{i_{1} \in \mathcal{S}}\sum_{\lambda_{i_{1}} \in \sum'_{i_{1}}} [\alpha^{2}_{i_{1}}(\lambda_{i_{1}}) + 2 \alpha_{i_{1}}(\lambda_{i_{1}})].rank(Gr^{i_{1}}_{\lambda_{i_{1}}}).[D_{i_{1}}]^{2}.[D]$
$$
$$
$\hspace{.6cm} - \frac{1}{2}\sum_{i_{1} \in \mathcal{S}}\sum_{\lambda_{i_{1}} \in \sum'_{i_{1}}} \alpha_{i_{1}}(\lambda_{i_{1}}).(\xi_{i_{1}})_{\star}\left((c_{1}^{D_{i_{1}}})^{2}(Gr^{i_{1}}_{\lambda_{i_{1}}})\right)$
$$
$$
$\hspace{.6cm} + \sum_{i_{1} \in \mathcal{S}}\sum_{\lambda_{i_{1}} \in \sum'_{i_{1}}} \alpha_{i_{1}}(\lambda_{i_{1}}).(\xi_{i_{1}})_{\star}\left(c_{2}^{D_{i_{1}}}(Gr_{\lambda_{i_{1}}}^{i_{1}})\right)$
$$
$$
$\hspace{.6cm} - \frac{1}{2} \sum_{i _{1} <i_{2}} \sum_{\aatop{\lambda_{i_{1}} }{ \lambda_{i_{2}}}} \sum_{p \in Irr(D_{i_{1}} \cap D_{i_{2}})}  [\alpha^{2}_{i_{2}}(\lambda_{i_{2}}).\alpha_{i_{1}}(\lambda_{i_{1}}) + 2 \alpha_{i_{1}}(\lambda_{i_{1}}).\alpha_{i_{2}}(\lambda_{i_{2}}) + \alpha_{i_{1}}(\lambda_{i_{1}}) + \alpha^{2}_{i_{2}}(\lambda_{i_{2}}) + 2\alpha_{i_{2}}(\lambda_{i_{2}})] .rank_{p}(Gr^{i_{1},i_{2}}_{\lambda_{i_{1}},\lambda_{i_{2}}}).[D_{i_{2}}].[D_{p}]$ 
$$
$$
$\hspace{.6cm} - \frac{1}{2} \sum_{i _{1} <i_{2}} \sum_{\aatop{\lambda_{i_{1}} }{ \lambda_{i_{2}}}}  \sum_{p \in Irr(D_{i_{1}} \cap D_{i_{2}})}  [\alpha^{2}_{i_{1}}(\lambda_{i_{1}}).\alpha_{i_{2}}(\lambda_{i_{2}}) + 2 \alpha_{i_{1}}(\lambda_{i_{1}}).\alpha_{i_{2}}(\lambda_{i_{2}}) + \alpha_{i_{2}}(\lambda_{i_{2}}) + \alpha^{2}_{i_{1}}(\lambda_{i_{1}}) + 2\alpha_{i_{1}}(\lambda_{i_{1}})] .rank_{p}(Gr^{i_{1},i_{2}}_{\lambda_{i_{1}},\lambda_{i_{2}}}).[D_{i_{1}}].[D_{p}]$ 
$$
$$
$\hspace{.6cm} + \sum_{i_{1} <i_{2}} \sum_{\aatop{\lambda_{i_{1}} }{ \lambda_{i_{2}}}}  \sum_{p \in Irr(D_{i_{1}} \cap D_{i_{2}})} [\alpha_{i_{1}}(\lambda_{i_{1}}).\alpha_{i_{2}}(\lambda_{i_{2}}) + \alpha_{i_{1}}(\lambda_{i_{1}}) + \alpha_{i_{2}}(\lambda_{i_{2}})].
rank_{p}(Gr^{i_{1},i_{2}}_{\lambda_{i_{1}},\lambda_{i_{2}}}).[D].[D_{p}]$ 
$$
$$
$ \hspace{.6cm} + \sum_{i_{1} <i_{2}} \sum_{\aatop{\lambda_{i_{1}}}{ \lambda_{i_{2}}}} [\alpha_{i_{1}}(\lambda_{i_{1}}).\alpha_{i_{2}}(\lambda_{i_{2}}) + 
\alpha_{i_{1}}(\lambda_{i_{1}}) + \alpha_{i_{2}}(\lambda_{i_{2}}) ].
(\xi_{i_{1},i_{2}})_{\star}\left(c_{1}^{D_{i_{1}} \cap D_{i_{2}}}(Gr^{i_{1},i_{2}}_{\lambda_{i_{1}},\lambda_{i_{2}}})\right)$
$$
$$
$ \hspace{.6cm} - \sum_{i_{1} <i_{2} < i_{3}} \sum_{\aatop{\lambda_{i_{1}} }{  \lambda_{i_{3}}}}  \sum_{p \in Irr(D_{i_{1}} \cap D_{i_{2}} \cap D_{i_{3}})}[\alpha_{i_{1}}(\lambda_{i_{1}}).\alpha_{i_{2}}(\lambda_{i_{2}}).\alpha_{i_{3}}(\lambda_{i_{3}}) + \alpha_{i_{1}}(\lambda_{i_{1}}).\alpha_{i_{2}}(\lambda_{i_{2}}) + $
\newline
$\alpha_{i_{2}}(\lambda_{i_{2}}).\alpha_{i_{3}}(\lambda_{i_{3}}) + \alpha_{i_{1}}(\lambda_{i_{1}}).\alpha_{i_{3}}(\lambda_{i_{3}}) + \alpha_{i_{1}}(\lambda_{i_{1}}) + \alpha_{i_{2}}(\lambda_{i_{2}}) + \alpha_{i_{3}}(\lambda_{i_{}})].rank_{p}(Gr^{i_{1},i_{2},i_{3}}_{\lambda_{i_{1}},\lambda_{i_{2}},\lambda_{i_{3}}}).[D_{p}].$

%%%%%%%%%%%%%%%%%%%%%%%%%%%%%%%%%%%%%%%%%%%%%%%%%%%%%%%%%%%%%%%%%%%%%%%%%%%%%%%%%%%%%%%%%%%%%%%%%%%%%%%%%%%%

\end{document}